\documentclass[a4, 11pt, leqno]{article}
\usepackage{amsmath}
\usepackage{amssymb}
\newtheorem{df}{Definition}[section]
\newtheorem{pr}[df]{Proposition}
\newtheorem{Th}[df]{Theorem}

\newtheorem{lm}[df]{Lemma}

\begin{document}

\title{\bf Solvability of an Initial-Boundary Value Problem for a Second Order Parabolic System with a Third Order Dispersion Term}
\author{Masashi A{\sc iki} and Tatsuo I{\sc guchi}}
\date{}
\maketitle
\vspace*{-0.5cm}
\begin{center}
Department of Mathematics, Faculty of Science and Technology, Keio University, \\
3-14-1 Hiyoshi, Kohoku-ku, Yokohama, 223-8522, J{\sc apan}
\end{center}

\begin{abstract}
We consider a linear second order parabolic system with a third order dispersion term. This type of system arises when considering a nonlinear model equation 
describing the motion of a vortex filament with axial flow immersed in an incompressible and inviscid fluid. 
We prove the solvability of an initial-boundary value problem of the parabolic-dispersive system
which allows application to the motion of a vortex filament. To do so, we propose a new regularization technique by adding a space-time derivative term.
\end{abstract}

\section{Introduction}
In this paper, we prove the unique solvability of the following initial-boundary value problems. For \( \alpha <0\),
\begin{eqnarray}
\left\{
\begin{array}{ll}
\mbox{\mathversion{bold}$u$}_{t}= \alpha \mbox{\mathversion{bold}$u$}_{xxx}+{\rm A}(\mbox{\mathversion{bold}$w$},\partial_{x})\mbox{\mathversion{bold}$u$} 
+ \mbox{\mathversion{bold}$f$}, &
x>0,t>0,\\
\mbox{\mathversion{bold}$u$}(x,0)=\mbox{\mathversion{bold}$u$}_{0}(x), & x>0, \\
\mbox{\mathversion{bold}$u$}_{x}(0,t)=\mbox{\mathversion{bold}$0$}, & t>0.
\end{array} \right.
\label{lpd1}
\end{eqnarray}
For \( \alpha >0\),
\begin{eqnarray}
\left\{
\begin{array}{ll}
\mbox{\mathversion{bold}$u$}_{t}= \alpha \mbox{\mathversion{bold}$u$}_{xxx}+{\rm A}(\mbox{\mathversion{bold}$w$},\partial_{x})\mbox{\mathversion{bold}$u$}
 + \mbox{\mathversion{bold}$f$}, &
x>0,t>0,\\
\mbox{\mathversion{bold}$u$}(x,0)=\mbox{\mathversion{bold}$u$}_{0}(x), & x>0, \\
\mbox{\mathversion{bold}$u$}(0,t)=\mbox{\mathversion{bold}$e$}, & t>0, \\
\mbox{\mathversion{bold}$u$}_{x}(0,t)=\mbox{\mathversion{bold}$0$}, & t>0.
\end{array} \right.
\label{lpd2}
\end{eqnarray}
Here, \( \mbox{\mathversion{bold}$u$}(x,t)=(u^{1}(x,t), u^{2}(x,t), \cdots  ,u^{m}(x,t)) \) is the unknown vector valued function, 
\( \mbox{\mathversion{bold}$u$}_{0}(x)\), \( \mbox{\mathversion{bold}$w$}(x,t)
= ( w^{1}(x,t),w^{2}(x,t), \cdots ,w^{k}(x,t))\), and \( \mbox{\mathversion{bold}$f$}(x,t)=( f^{1}(x,t), f^{2}(x,t), \cdots ,f^{m}(x,t))\) are known vector valued 
functions, \( \mbox{\mathversion{bold}$e$}\) is an arbitrary constant vector, subscripts denote derivatives with the respective variables, 
\( {\rm A}(\mbox{\mathversion{bold}$w$},\partial_{x})\) is 
a second order differential operator of the form \( {\rm A}({\mbox{\mathversion{bold}$w$}},\partial_{x})={\rm A}_{0}({\mbox{\mathversion{bold}$w$}})\partial ^{2}_{x} 
+ {\rm A}_{1}({\mbox{\mathversion{bold}$w$}})\partial _{x} + {\rm A}_{2}({\mbox{\mathversion{bold}$w$}}) \). \( {\rm A}_{0}, \ {\rm A}_{1}, \ {\rm A}_{2}\) are
smooth matrices and \( {\rm A}({\mbox{\mathversion{bold}$w$}},\partial_{x}) \) is strongly elliptic in the sense that for any bounded domain \( E\) in 
\( \mathbf{R}^{k}\), there is a positive constant \( \delta \) such that 
\begin{eqnarray*}
\inf\limits_{\mbox{\mathversion{bold}\scriptsize $w$}\in E}\left\{{\rm A}_{0}(\mbox{\mathversion{bold}$w$})
+{\rm A}_{0}(\mbox{\mathversion{bold}$w$})^{*}\right\}\geq \delta {\rm I},
\end{eqnarray*}
where \( {\rm I}\) is the unit matrix and \( *\) denotes the adjoint of a matrix.
The equation is a second order strongly parabolic equation with a third order constant coefficient dispersion term. This type of equation comes up 
when analyzing a model equation describing the motion of a vortex filament with axial flow. A vortex filament is a space curve on which the vorticity of an
incompressible and inviscid fluid is concentrated.
The authors are especially interested in the following initial-boundary value problems for the motion of a vortex filament. 
\begin{eqnarray}
\left\{
\begin{array}{ll}
\mbox{\mathversion{bold}$v$}_{t}=\mbox{\mathversion{bold}$v$}\times \mbox{\mathversion{bold}$v$}_{ss} 
+ \alpha \left\{ \mbox{\mathversion{bold}$v$}_{sss} + \frac{3}{2}\mbox{\mathversion{bold}$v$}_{ss}\times (\mbox{\mathversion{bold}$v$}\times \mbox{\mathversion{bold}$v$}_{s})
+ \frac{3}{2}\mbox{\mathversion{bold}$v$}_{s}\times ( \mbox{\mathversion{bold}$v$}\times \mbox{\mathversion{bold}$v$}_{ss}) \right\}, & s>0,t>0,\\
\mbox{\mathversion{bold}$v$}(s,0)=\mbox{\mathversion{bold}$v$}_{0}(s),& s>0, \\
\mbox{\mathversion{bold}$v$}_{s}(0,t)= \mbox{\mathversion{bold}$0$},& t>0,
\end{array}\right.
\label{vf1}
\end{eqnarray}
\begin{eqnarray}
\left\{
\begin{array}{ll}
\mbox{\mathversion{bold}$v$}_{t}=\mbox{\mathversion{bold}$v$}\times \mbox{\mathversion{bold}$v$}_{ss} 
+ \alpha \left\{ \mbox{\mathversion{bold}$v$}_{sss} + \frac{3}{2}\mbox{\mathversion{bold}$v$}_{ss}\times (\mbox{\mathversion{bold}$v$}\times \mbox{\mathversion{bold}$v$}_{s})
+ \frac{3}{2}\mbox{\mathversion{bold}$v$}_{s}\times ( \mbox{\mathversion{bold}$v$}\times \mbox{\mathversion{bold}$v$}_{ss}) \right\}, & s>0,t>0,\\
\mbox{\mathversion{bold}$v$}(s,0)=\mbox{\mathversion{bold}$v$}_{0}(s),& s>0, \\
\mbox{\mathversion{bold}$v$}(0,t)=\mbox{\mathversion{bold}$e$}_{3},& t>0,\\
\mbox{\mathversion{bold}$v$}_{s}(0,t)= \mbox{\mathversion{bold}$0$},& t>0,
\end{array}\right.
\label{vf2}
\end{eqnarray}
where \( \mbox{\mathversion{bold}$v$}=(v^{1}(s,t),v^{2}(s,t),v^{3}(s,t))\) is the tangent vector of the filament parameterized by its arc length \( s\) at time \( t\), 
\( \mbox{\mathversion{bold}$e$}_{3}\) is the unit upward vector, \( \times \) is the exterior product in the three dimensional Euclidean space, and \( \alpha \) is a real constant
describing the magnitude of the effect of axial flow. As far as the authors know, there are no results on initial-boundary value problems for the above equation. 
In Nishiyama and Tani \cite{5}, they proved the unique solvability globally in time of the Cauchy problem in Sobolev spaces. Onodera \cite{6,7} proved the unique solvability 
of the Cauchy problem for a geometrically generalized equation. Segata \cite{8} proved the unique solvability and showed the asymptotic behavior in time of the solution to the 
Hirota equation, given by
\begin{eqnarray*}
{\rm i}q_{t}=q_{xx}+\frac{1}{2}|q|^{2}q + {\rm i}\alpha \big( -q_{xxx} + |q|^{2}q_{x}\big),
\end{eqnarray*}
which can be obtained by applying the Hasimoto transformation to the vortex filament equation. 
Although there are many literatures regarding Schr\"odinger type equations, for (\ref{vf2}), the boundary condition does not transfer into 
a form that is manageable, so we decided to work with the vortex filament equation directly.

Using the fact that a smooth solution satisfies \( |\mbox{\mathversion{bold}$v$}|=1\), linearizing around \( \mbox{\mathversion{bold}$w$}\), and
neglecting lower order terms, we obtain
\begin{eqnarray*}
\mbox{\mathversion{bold}$v$}_{t}= \mbox{\mathversion{bold}$w$}\times \mbox{\mathversion{bold}$v$}_{ss} + \alpha \big\{ \mbox{\mathversion{bold}$v$}_{sss}
+3\mbox{\mathversion{bold}$v$}_{ss}\times (\mbox{\mathversion{bold}$w$}\times \mbox{\mathversion{bold}$w$}_{s})\big\} + \mbox{\mathversion{bold}$f$}.
\end{eqnarray*}
We note that the second order derivative terms have skew-symmetric coefficients, so if we regularize the equation with a second order viscosity term, we have a parabolic-dispersive system
given by
\begin{eqnarray}
\mbox{\mathversion{bold}$v$}_{t}= \alpha \mbox{\mathversion{bold}$v$}_{sss} + \big\{\delta \mbox{\mathversion{bold}$v$}_{ss}
+ \mbox{\mathversion{bold}$w$}\times \mbox{\mathversion{bold}$v$}_{ss} +
3\alpha \mbox{\mathversion{bold}$v$}_{ss}\times (\mbox{\mathversion{bold}$w$}\times \mbox{\mathversion{bold}$w$}_{s})\big\} + \mbox{\mathversion{bold}$f$},
\label{vf3}
\end{eqnarray}
and the solvability of the initial-boundary value problem for the above system plays a crucial role in solving (\ref{vf1}) and (\ref{vf2}).
This is our motivation for considering (\ref{lpd1}) and (\ref{lpd2}). The application of the results of this paper to the 
vortex filament equation will be considered in a forth coming paper.

At first glance, one may think that (\ref{vf3}) can be treated by using known results of KdV and KdV-Burgers equations such as 
Hayashi and Kaikina \cite{3}, Hayashi, Kaikina, and Ruiz Paredes \cite{4}, or Bona and Zhang \cite{9}. This seems hard to do because the vortex filament equation 
in (\ref{vf1}) and (\ref{vf2}) has 
second order derivatives in the nonlinear term and the linear estimates obtained in the KdV and KdV-Burgers theory 
is not enough to treat the nonlinear term as a regular perturbation. Thus, a need for a new technique arose.

Our method to prove the solvability of (\ref{lpd1}) and (\ref{lpd2}) is parabolic regularization. For (\ref{lpd2}), we can regularize the equation with a 
fourth order viscosity term, making it a standard parabolic system. 
We can not do this for (\ref{lpd1}) because a fourth order parabolic equation requires two boundary conditions to solve, but 
our original problem imposes only one boundary condition. Thus, a standard regularization can not be applied to (\ref{lpd1}). 
To prove the unique solvability of (\ref{lpd1}), we introduce a new type of regularization, namely, we consider the following regularized problem.
\begin{eqnarray}
\left\{
\begin{array}{ll}
\mbox{\mathversion{bold}$u$}_{t}= \alpha \left(  \mbox{\mathversion{bold}$u$}_{xx} - \varepsilon \mbox{\mathversion{bold}$u$}_{t}\right)_{x}
+ {\rm A}(\mbox{\mathversion{bold}$w$},\partial _{x})\mbox{\mathversion{bold}$u$}+ \mbox{\mathversion{bold}$g$} ,& x>0, t>0, \\
\mbox{\mathversion{bold}$u$}(x,0)= \mbox{\mathversion{bold}$u$}_{0}(x), & x>0, \\
\mbox{\mathversion{bold}$u$}_{x}(0,t)=\mbox{\mathversion{bold}$0$}, & t>0,
\end{array} \right.
\label{ep}
\end{eqnarray}
where \( \varepsilon >0\). To construct the solution of the above system, we first consider the following problem.
\begin{eqnarray}
\left\{
\begin{array}{ll}
\mbox{\mathversion{bold}$u$}_{t}= \alpha \left(  \mbox{\mathversion{bold}$u$}_{xx} - \varepsilon \mbox{\mathversion{bold}$u$}_{t}\right)_{x}
+ \mbox{\mathversion{bold}$g$} ,& x>0, t>0, \\
\mbox{\mathversion{bold}$u$}(x,0)= \mbox{\mathversion{bold}$u$}_{0}(x), & x>0, \\
\mbox{\mathversion{bold}$u$}_{x}(0,t)=\mbox{\mathversion{bold}$0$}, & t>0,
\end{array} \right.
\label{lp1}
\end{eqnarray}
(\ref{ep}) is a parabolic regularization of (\ref{lpd1}) and the principle terms are the terms in parenthesis. In fact if we substitute \( \mbox{\mathversion{bold}$u$}(x,t)=
{\rm e}^{\tau t + {\rm i}\xi x}\) into \( \mbox{\mathversion{bold}$u$}_{t}= 
\alpha \left(  \mbox{\mathversion{bold}$u$}_{xx} - \varepsilon \mbox{\mathversion{bold}$u$}_{t}\right)_{x}\), we obtain the dispersion relation \( \tau = -\alpha (\xi ^{2} 
+ \varepsilon \tau ){\rm i}\xi \), so that
for a non-trivial solution to exist, we need
\begin{eqnarray*}
\Re{\tau } = -\frac{\alpha ^{2}\varepsilon \xi ^{4}}{1+\alpha ^{2}\varepsilon ^{2}\xi ^{2}},
\end{eqnarray*}
which indicates that the equation is parabolic in nature.

Because the proof for the case \( \alpha >0\) is fairly standard, we concentrate on the case \( \alpha <0\), and give a remark on the case \( \alpha >0\)
at the end.

The contents of this paper are as follows. In section 2, we define function spaces and notations that are used in this paper. In section 3, we look at the
compatibility conditions and the necessary corrections to the given data required for the regularized system. Then in section 4, we construct and estimate 
the solution to the regularized system. In section 5, we construct and estimate the solution of the parabolic-dispersive system in appropriate function spaces
and give our main theorem of this paper. 
Finally in section 6, we give a remark on the case \( \alpha >0\).


\section{Function Spaces}

We define some function spaces that will be used throughout this paper, and notations associated with the spaces.
For an open interval \( \Omega \), 
a non-negative integer \( m\), and \( 1\leq p \leq \infty \), \( W^{m,p}(\Omega )\) is the Sobolev space 
containing all real-valued functions that have derivatives in the sense of distribution up to order \( m\) 
belonging to \( L^{p}(\Omega )\) and \( \dot{W}^{m,p}(\Omega ) \) is the homogeneous Sobolev space.
We set \( H^{m}(\Omega ) := W^{m,2}(\Omega ) \) as the Sobolev space equipped with the usual inner product 
and \( \dot{H}^{m}(\Omega ):=\dot{W}^{m,2}(\Omega ) \).  
We will particularly use the cases \( \Omega= {\mathbf R} \) and \( \Omega = {\mathbf R}_{+} \), 
where \( {\mathbf R}_{+} = \{ x\in {\mathbf R}; x>0 \} \). 
When \(\Omega = \mathbf{R}_{+}\), the norm in \( H^{m}(\Omega ) \) is denoted by \( \| \cdot \|_{m} \) and we simply write \( \| \cdot \| \) 
for \( \|\cdot \|_{0} \). Otherwise, for a Banach space \( X\), the norm in \( X\) is written as \( \| \cdot \| _{X}\).
The inner product in \( L^{2}(\mathbf{R}_{+})\) is denoted by \( (\cdot ,\cdot )\) 
and the inner product in \( L^{2}(\mathbf{R})\) is denoted by \( \langle\cdot, \cdot \rangle\).

For \( 0<T< \infty \) and a Banach space \( X\), 
\( C^{m}([0,T];X) \) denotes the space of functions that are \( m\) times continuously differentiable 
in \( t\) with respect to the norm of \( X\).

We define the Sobolev--Slobodetski\u\i \ space. For \( 0<T\leq \infty \), we denote 
\( Q_{T}:= \mathbf{R}_{+}\times (0,T)\) and for \( h>0\) and a positive integer \( l \), we define the space 
\( H^{l,l/2}_{h}(Q_{T})\) as the space of functions defined on \( Q_{T}\) with finite norm 
\begin{eqnarray*}
||| u||| ^{2}_{H^{l,l/2}_{h}(Q_{T})} := ||| u||| ^{2}_{H^{l,0}_{h}(Q_{T})} + ||| u||| ^{2}_{H^{0,l/2}_{h}(Q_{T})},
\end{eqnarray*}
where
\begin{eqnarray*}
\begin{aligned}
||| u||| ^{2}_{H^{l,0}_{h}(Q_{T})}&:= \int^{T}_{0}{\rm e}^{-2ht}\| u(\cdot ,t)\| ^{2}_{\dot{H}^{l}}{\rm d}t, \\
||| u||| ^{2}_{H^{0,l/2}_{h}(Q_{T})}&:= h^{l}\int^{T}_{0}{\rm e}^{-2ht}\| u(\cdot ,t)\|^{2}{\rm d}t \\
&\hspace*{1cm}+\int^{T}_{0}{\rm e}^{-2ht}\int^{\infty}_{0}\left\| \frac{\partial^{[l/2]}u_{0}(\cdot ,t-r )}{\partial t^{[l/2]}}
           -\frac{\partial^{[l/2]}u_{0}(\cdot ,t)}{\partial t^{[l/2]}}\right\| ^{2}r ^{-1-l+2[\frac{l}{2}]}{\rm d}r {\rm d}t,
\end{aligned}
\end{eqnarray*}
\( [\frac{l}{2}]\) is the integer part of \( \frac{l}{2}\) and \( u_{0}\) is the extension of \( u\) by zero into \( t<0\) if \( \frac{l}{2}\) is not an integer. 
When \( \frac{l}{2}\) is an integer,
\begin{eqnarray*}
||| u||| ^{2}_{H^{0,l/2}_{h}(Q_{T})}:= \int^{T}_{0}{\rm e}^{-2ht}\left( h^{l}\| u(\cdot ,t)\| ^{2}+\left\| \frac{\partial^{l/2}u}{\partial t^{l/2}}(\cdot ,t)\right\| ^{2}\right) {\rm d}t
\end{eqnarray*}
and we also impose that \( \frac{\partial ^{j}u}{\partial t^{j}}(x,0)=0\) for \( j=0,1, \dots , \frac{l}{2}-1\). 
When \( T=\infty \), the following equivalent norm for the space \( H^{l,l/2}_{h}(Q_{\infty})\) will be used.
\begin{eqnarray*}
\| u\|  ^{2}_{H^{l,l/2}_{h}(Q_{\infty})}:= \sum _{j\leq l} \int^{\infty}_{-\infty}\left \| 
\frac{\partial ^{j}\tilde{u}}{\partial x^{j}}(\cdot ,\tau ) \right\|^{2}|\tau |^{l-j}{\rm d}\eta ,
\end{eqnarray*}
where a tilde denotes the Laplace transform in \( t\) defined by 
\begin{eqnarray*}
\tilde{u}(x,\tau ) = \int ^{\infty}_{0} {\rm e}^{-\tau t}u(x,t){\rm d}t,
\end{eqnarray*}
where \( \tau =h+{\rm i}\eta \) with \( h>0\). The equivalence is shown in Solonnikov \cite{2}.

For any function space described above, we say that a vector valued function belongs to the function space 
if each of its components does.

%
%
\section{Compatibility Conditions}
\setcounter{equation}{0}
We will construct the solution of (\ref{lpd1}) by taking the limit \( \varepsilon \rightarrow 0\) in the following regularized system.
\begin{eqnarray}
\left\{
\begin{array}{ll}
\mbox{\mathversion{bold}$u$}_{t}= -\alpha \varepsilon \mbox{\mathversion{bold}$u$}_{tx} + 
\alpha \mbox{\mathversion{bold}$u$}_{xxx}+{\rm A}(\mbox{\mathversion{bold}$w$},\partial_{x})\mbox{\mathversion{bold}$u$} + \mbox{\mathversion{bold}$g$}, &
x>0,t>0,\\
\mbox{\mathversion{bold}$u$}(x,0)=\mbox{\mathversion{bold}$u$}_{0}(x), & x>0, \\
\mbox{\mathversion{bold}$u$}_{x}(0,t)=\mbox{\mathversion{bold}$0$}, & t>0.
\end{array} \right.
\label{lpp1}
\end{eqnarray}
Since the derivation of the compatibility conditions for the regularized system is complicated and the required corrections for the given data is not standard, 
we devote this section to clarify these matters.

\subsection{Compatibility Conditions for (\ref{lpd1})}
We first define the compatibility condition for the original system (\ref{lpd1}). We denote the equation in (\ref{lpd1}) as
\begin{eqnarray}
\mbox{\mathversion{bold}$Q$}_{1}(\mbox{\mathversion{bold}$u$},\mbox{\mathversion{bold}$f$},\mbox{\mathversion{bold}$w$}) = \alpha \mbox{\mathversion{bold}$u$}_{xxx}
+{\rm A}(\mbox{\mathversion{bold}$w$},\partial _{x})\mbox{\mathversion{bold}$u$} + \mbox{\mathversion{bold}$f$} ,
\label{ccq1}
\end{eqnarray}
and we also use the notation 
\( \mbox{\mathversion{bold}$Q$}_{1}(x,t):=\mbox{\mathversion{bold}$Q$}_{1}(\mbox{\mathversion{bold}$u$},\mbox{\mathversion{bold}$f$},\mbox{\mathversion{bold}$w$}) \)
and sometimes we omit the \( (x,t)\) for simplicity. We successively define
\begin{eqnarray}
\mbox{\mathversion{bold}$Q$}_{n}:= \alpha \partial^{3}_{x}\mbox{\mathversion{bold}$Q$}_{n-1}
+ \sum^{n-1}_{j=0}\left(
\begin{array}{c}
n-1\\
j
\end{array}\right)
{\rm B}_{j}
\mbox{\mathversion{bold}$Q$}_{n-1-j} + \partial ^{n-1}_{t}\mbox{\mathversion{bold}$f$},
\label{ccqn}
\end{eqnarray}
where \( {\rm B}_{j}= \big( \partial ^{j}_{t}{\rm A_{0}}(\mbox{\mathversion{bold}$w$}) \big) \partial^{2}_{x} + \big( \partial ^{j}_{t}{\rm A_{1}}(\mbox{\mathversion{bold}$w$})\big)
\partial_{x} + \partial^{j}_{t}{\rm A_{2}}(\mbox{\mathversion{bold}$w$})\).
The above definition gives the formula for the expression of \( \partial ^{n}_{t}\mbox{\mathversion{bold}$u$}\) 
which only contains \( x\) derivatives of \( \mbox{\mathversion{bold}$u$}\) and mixed derivatives of \( \mbox{\mathversion{bold}$w$}\) and 
\( \mbox{\mathversion{bold}$f$}\). From the boundary condition in (\ref{lpd1}), we arrive at the following 
definition for the compatibility conditions.
\begin{df}{\rm (}Compatibility conditions for {\rm (\ref{lpd1})}{\rm )}.
For \( n\in \mathbf{N}\cup \{ 0\}\), we say that \( \mbox{\mathversion{bold}$u$}_{0}\), \( \mbox{\mathversion{bold}$f$}\), and \( \mbox{\mathversion{bold}$w$}\) satisfy the 
\(n\)-th order compatibility condition for {\rm (\ref{lpd1})} if
\begin{eqnarray*}
\mbox{\mathversion{bold}$u$}_{0x}(0,0)=\mbox{\mathversion{bold}$0$}
\end{eqnarray*}
when \( n=0\), and
\begin{eqnarray*}
\big( \partial _{x}\mbox{\mathversion{bold}$Q$}_{n}\big)(0,0)=\mbox{\mathversion{bold}$0$}
\end{eqnarray*}
when \( n\geq 1\).
We also say that the data satisfy the compatibility conditions for {\rm (\ref{lpd1})} up to order \( n\) 
if they satisfy the \( k\)-th order compatibility condition for all \( k\) with \( 0\leq k\leq n\).
\label{cclpd1}
\end{df}

Now that we have defined the compatibility conditions, we discuss an approximation of the data via smooth functions which keep the compatibility conditions.
The function spaces we consider for the data and solution of (\ref{lpd1}) are the following. Let \( l\) be a non-negative integer. 
\( X^{l}\) is the function space that we are constructing the solution in, specifically,
\begin{eqnarray*}
X^{l}:= \bigcap ^{l}_{j=0} \bigg( C^{j}\big( [0,T];H^{2+3(l-j)}(\mathbf{R}_{+})\big) \cap H^{j}\big( 0,T;H^{3+3(l-j)}(\mathbf{R}_{+})\big) \bigg).
\end{eqnarray*}
As a consequence, \( \mbox{\mathversion{bold}$u$}_{0}\) will be required to belong in \( H^{2+3l}(\mathbf{R}_{+})\). 
\( Y^{l}\) is the function space that \( {\mbox{\mathversion{bold}$f$}}\) will be required to belong in, and is defined by 
\begin{eqnarray*}
Y^{l}:= \bigg\{ f; \  f\in \bigcap ^{l-1}_{j=0} C^{j}\big( [0,T];H^{2+3(l-1-j)}(\mathbf{R}_{+})\big) , \ \partial ^{l}_{t}f \in L^{2}\big( 0,T; H^{1}(\mathbf{R}_{+})\big) \bigg\}.
\end{eqnarray*}
Finally, \( Z^{l}\) is the function space that \( {\mbox{\mathversion{bold}$w$}}\) will belong in and is defined as
\begin{eqnarray*}
Z^{l}:= \bigg\{ w; \ w\in \bigcap ^{l-1}_{j=0}C^{j}\big( [0,T];H^{2+3(l-1-j)}(\mathbf{R}_{+})\big), \ \partial ^{l}_{t}w\in L^{\infty}\big (0,T;H^{1}(\mathbf{R}_{+})\big) \bigg\}.
\end{eqnarray*}
Data belonging to the above function spaces with index \( l\) are smooth enough for the \( l\)-th order compatibility condition to have meaning in a point-wise sense,
but the \( (l+1)\)-th order compatibility condition does not.
By utilizing the method in Rauch and Massey \cite{1}, we can get the following.
\begin{lm}
Fix non-negative integers \( l\) and \( N\) with \( N>l\). 
For any \( \mbox{\mathversion{bold}$u$}_{0}\in H^{2+3l}(\mathbf{R}_{+})\), \( \mbox{\mathversion{bold}$f$}\in Y^{l}\), and 
\( \mbox{\mathversion{bold}$w$}\in Z^{l}\) satisfying the compatibility conditions for {\rm (\ref{lpd1})} up to order \( l\), there exist sequences 
\( \{ \mbox{\mathversion{bold}$u$}_{0n}\}_{n\geq1}\subset H^{2+3N}(\mathbf{R}_{+}) \),
\( \{ \mbox{\mathversion{bold}$f$}_{n}\}_{n\geq 1}\subset Y^{N}\), and \( \{\mbox{\mathversion{bold}$w$}_{n}\}_{n\geq 1}\subset Z^{N}\) such that 
for any \( n\geq 1\), \( \mbox{\mathversion{bold}$u$}_{0n}\), \( \mbox{\mathversion{bold}$f$}_{n}\), and \( \mbox{\mathversion{bold}$w$}_{n}\) satisfy 
the compatibility conditions for {\rm (\ref{lpd1})} up to order \( N\) and 
\begin{eqnarray*}
\mbox{\mathversion{bold}$u$}_{0n}\rightarrow \mbox{\mathversion{bold}$u$}_{0} \ in \ H^{2+3l}(\mathbf{R}_{+}), \ \ 
\mbox{\mathversion{bold}$f$}_{n}\rightarrow \mbox{\mathversion{bold}$f$} \ in \ Y^{l}, \ and \ 
\mbox{\mathversion{bold}$w$}_{n}\rightarrow \mbox{\mathversion{bold}$w$} \ in \ Z^{l}.
\end{eqnarray*}
\label{appr}
\end{lm}
From Lemma \ref{appr}, we can assume that the given data are arbitrarily smooth and satisfy the necessary compatibility conditions in the
proceeding arguments.

\subsection{Compatibility Conditions for (\ref{lpp1})}
In this subsection, we define the compatibility conditions for (\ref{lpp1}). We write the equation in (\ref{lpp1}) as
\begin{eqnarray}
\mbox{\mathversion{bold}$u$}_{t}=- \alpha \varepsilon \mbox{\mathversion{bold}$u$}_{tx} 
+ \mbox{\mathversion{bold}$P$}_{1}(\mbox{\mathversion{bold}$u$},\mbox{\mathversion{bold}$g$},\mbox{\mathversion{bold}$w$}),
\label{ccpp1}
\end{eqnarray}
in other words, \( \mbox{\mathversion{bold}$P$}_{1}(\mbox{\mathversion{bold}$u$}, \mbox{\mathversion{bold}$g$}, \mbox{\mathversion{bold}$w$})=
\alpha \mbox{\mathversion{bold}$u$}_{xxx} + {\rm A}(\mbox{\mathversion{bold}$w$},\partial _{x})\mbox{\mathversion{bold}$u$} + \mbox{\mathversion{bold}$g$}\). 
We use the notations \( \mbox{\mathversion{bold}$P$}_{1}(x,t)\) and \( \mbox{\mathversion{bold}$P$}_{1}\) as we did with \( \mbox{\mathversion{bold}$Q$}_{1}\) in the 
last subsection.
Setting \( \mbox{\mathversion{bold}$\phi$} _{1} (x):= \mbox{\mathversion{bold}$u$}_{t}(x,0)\) and taking the trace \( t=0\) of the equation we have
\begin{eqnarray}
\alpha \varepsilon \mbox{\mathversion{bold}$\phi$} _{1} ^{\prime} + \mbox{\mathversion{bold}$\phi$} _{1} 
=\mbox{\mathversion{bold}$P$}_{1}(\cdot ,0).
\label{cc1}
\end{eqnarray}
A prime denotes a derivative with respect to \( x\). 
Note that \( \mbox{\mathversion{bold}$P$}_{1}(x,0)\) is expressed using given data only.
Solving the above ordinary differential equation for \( \mbox{\mathversion{bold}$\phi$} _{1}\) we have
\begin{eqnarray*}
\mbox{\mathversion{bold}$\phi$} _{1}(x)= {\rm e}^{-\frac{x}{\alpha \varepsilon}}\left\{ \mbox{\mathversion{bold}$\phi$} _{1}(0)+\frac{1}{\alpha \varepsilon} 
\int^{x}_{0}{\rm e}^{\frac{y}{\alpha \varepsilon}}
\mbox{\mathversion{bold}$P$}_{1}(y,0) {\rm d}y\right\}. 
\end{eqnarray*}
Since we are looking for solutions that are square integrable, we impose that $ \displaystyle\lim _{x \to \infty}  \mbox{\mathversion{bold}$\phi$} _{1}(x) =
 \mbox{\mathversion{bold}$0$} $, so we have
\begin{eqnarray*}
\mbox{\mathversion{bold}$\phi$} _{1}(0)=-\frac{1}{\alpha \varepsilon} \int^{\infty}_{0}{\rm e}^{\frac{y}{\alpha \varepsilon}}
{\mbox{\mathversion{bold}$P$}}_{1}(y,0){\rm d}y,
\end{eqnarray*}
which gives
\begin{eqnarray*}
\mbox{\mathversion{bold}$\phi$} _{1}(x) = -\frac{1}{\alpha \varepsilon}\int ^{\infty}_{x}{\rm e}^{-\frac{1}{\alpha \varepsilon}(x-y)}
{\mbox{\mathversion{bold}$P$}}_{1}(y,0) {\rm d}y.
\end{eqnarray*}
By direct calculation, we see that 
\begin{eqnarray*}
\mbox{\mathversion{bold}$\phi$} _{1}^{\prime}(x) = -\frac{1}{\alpha \varepsilon}\int ^{\infty}_{x}{\rm e}^{-\frac{1}{\alpha \varepsilon}(x-y)}
{\mbox{\mathversion{bold}$P$}}_{1}^{\prime}(y,0) {\rm d}y,
\end{eqnarray*}
where we have used integration by parts. We also note here that \( \mbox{\mathversion{bold}$\phi$}_{1}\) is expressed with given data only. 
From the boundary condition in (\ref{lpp1}), we see that the first order compatibility condition is
\begin{eqnarray*}
\int ^{\infty}_{0}{\rm e}^{\frac{y}{\alpha \varepsilon}}
{\mbox{\mathversion{bold}$P$}}_{1}^{\prime}(y,0) {\rm d}y = \mbox{\mathversion{bold}$0$}.
\end{eqnarray*}
In the same manner, we will derive the \( n\)-th order compatibility condition for \( n\geq 2\).
Taking the \( t\) derivative of the equation in (\ref{lpp1}) \( (n-1)\) times, 
taking the trace \( t=0\), and setting \( \mbox{\mathversion{bold}$\phi$}_{n}(x):= \partial ^{n}_{t}\mbox{\mathversion{bold}$u$}(x,0)\), we have
\begin{eqnarray*}
\alpha \varepsilon \mbox{\mathversion{bold}$\phi$}_{n}^{\prime}+ \mbox{\mathversion{bold}$\phi $}_{n}= \partial ^{n-1}_{t}\mbox{\mathversion{bold}$P$}_{1}.
\end{eqnarray*}
We denote 
\begin{eqnarray*}
\mbox{\mathversion{bold}$P$}_{n}:= \partial ^{n-1}_{t}\mbox{\mathversion{bold}$P$}_{1}.
\end{eqnarray*}
We will prove by induction that \( \mbox{\mathversion{bold}$\phi$}_{n}\) and \( \mbox{\mathversion{bold}$P$}_{n}(x,0)\) are expressed using given data only.
Since \( \mbox{\mathversion{bold}$P$}_{n}= \partial ^{n-1}_{t}\mbox{\mathversion{bold}$P$}_{n-1} = \partial ^{n-1}_{t}( \alpha \mbox{\mathversion{bold}$u$}_{xxx}
+ {\rm A}(\mbox{\mathversion{bold}$w$})\mbox{\mathversion{bold}$u$} + \mbox{\mathversion{bold}$g$})\), it holds that
\begin{eqnarray}
\mbox{\mathversion{bold}$P$}_{n}(\cdot,0) = \alpha \mbox{\mathversion{bold}$\phi$}_{n-1}^{\prime \prime \prime}+ \sum ^{n-1}_{j=0}
\left(
\begin{array}{c}
n-1\\
j
\end{array}\right)
{\rm B}_{j}\mbox{\mathversion{bold}$\phi$}_{n-1-j} + \partial ^{n-1}_{t}\mbox{\mathversion{bold}$g$}(\cdot, 0).
\label{ccppn}
\end{eqnarray}
For a \( n\geq 2\), assume that \( \mbox{\mathversion{bold}$\phi$}_{k}\) and \( \mbox{\mathversion{bold}$P$}_{k}(x,0)\) are expressed with given data
for \( 1\leq k\leq n-1\).
Solving for \( \mbox{\mathversion{bold}$\phi $}_{n}\) yields
\begin{eqnarray*}
\mbox{\mathversion{bold}$\phi$} _{n}(x) = -\frac{1}{\alpha \varepsilon}\int ^{\infty}_{x}{\rm e}^{-\frac{1}{\alpha \varepsilon}(x-y)}
{\mbox{\mathversion{bold}$P$}}_{n}(y,0) {\rm d}y.
\end{eqnarray*}
This proves that \( \mbox{\mathversion{bold}$P$}_{n}(x,0)\) and \( \mbox{\mathversion{bold}$\phi$}_{n}\) are 
expressed using given data only.

Again by direct calculation, we have
\begin{eqnarray*}
\mbox{\mathversion{bold}$\phi$} _{n}^{\prime}(x) = -\frac{1}{\alpha \varepsilon}\int ^{\infty}_{x}{\rm e}^{-\frac{1}{\alpha \varepsilon}(x-y)}
{\mbox{\mathversion{bold}$P$}}_{n}^{\prime}(y,0) {\rm d}y,
\end{eqnarray*}
and arrive at the \( n\)-th order compatibility condition 
\begin{eqnarray*}
\int ^{\infty}_{0}{\rm e}^{\frac{y}{\alpha \varepsilon}}
{\mbox{\mathversion{bold}$P$}}_{n}^{\prime}(y,0) {\rm d}y = \mbox{\mathversion{bold}$0$}.
\end{eqnarray*}
Now we can define the following.
\begin{df}
{\rm (}Comaptibility conditions for {\rm(\ref{lpp1}))}. For \( n\in \mathbf{N}\cup\{0\}\), we say that 
\( \mbox{\mathversion{bold}$u$}_{0}\), \( \mbox{\mathversion{bold}$g$}\), and \( \mbox{\mathversion{bold}$w$}\) satisfy 
the \( n\)-th order compatibility condition for {\rm (\ref{lpp1})} if 
\begin{eqnarray*}
\mbox{\mathversion{bold}$u$}_{0x}(0)=\mbox{\mathversion{bold}$0$}
\end{eqnarray*}
when \( n=0\), and 
\begin{eqnarray*}
\int ^{\infty}_{0}{\rm e}^{\frac{y}{\alpha \varepsilon}}
{\mbox{\mathversion{bold}$P$}}_{n}^{\prime}(y,0) {\rm d}y = \mbox{\mathversion{bold}$0$}
\end{eqnarray*}
when \( n\geq 1\). We also say that the data satisfy the compatibility conditions for {\rm (\ref{lpp1})} up to order \( n\) if
the data satisfy the \( k\)-th order compatibility condition for all \( k\) with \( 0\leq k\leq n\).
For the definition of \( {\mbox{\mathversion{bold}$P$}}_{n}\), see {\rm (\ref{ccpp1})} and {\rm (\ref{ccppn})}.
\label{cclpp1}
\end{df}
We note that for \( \mbox{\mathversion{bold}$u$}_{0}\in H^{2+3l}(\mathbf{R}_{+})\), \( \mbox{\mathversion{bold}$f$}\in Y^{l}\), and
\( \mbox{\mathversion{bold}$w$}\in Z^{l}\), the compatibility conditions up to order \( l\) have
meaning in the point-wise sense, but the \( (l+1)\)-th order compatibility condition does not.

\subsection{Corrections to the Data}

Since we regularized the equation, we must make corrections to the data to assure that the compatibility conditions continue to hold. 
Fix a large integer \( N\) and suppose that \( {\mbox{\mathversion{bold}$u$}}_{0}\in H^{2+3N}(\mathbf{R}_{+})\), \( {\mbox{\mathversion{bold}$f$}}\in Y^{N}\),
and \( \mbox{\mathversion{bold}$w$}\in Z^{N}\) satisfy the compatibility conditions for (\ref{lpd1}) up to order \( N\). 
We will make corrections to the forcing term so that the data satisfy the compatibility conditions for (\ref{lpp1}) up to order \( N\). 
More specifically, we prove the following
\begin{pr}
Fix a positive integer \( N\). For \( {\mbox{\mathversion{bold}$u$}}_{0}\in H^{2+3N}(\mathbf{R}_{+})\), \( {\mbox{\mathversion{bold}$f$}}\in Y^{N}\),
and \( \mbox{\mathversion{bold}$w$}\in Z^{N}\) satisfying the compatibility conditions for {\rm (\ref{lpd1})} up to order \( N\), 
we can define \( \mbox{\mathversion{bold}$g$}\in Y^{N}\) in the form \( \mbox{\mathversion{bold}$g$}=\mbox{\mathversion{bold}$f$}+\mbox{\mathversion{bold}$h$}_{\varepsilon }\)
such that \( \mbox{\mathversion{bold}$u$}_{0}\), \( \mbox{\mathversion{bold}$g$}\), and \( \mbox{\mathversion{bold}$w$}\) satisfy the compatibility conditions 
for {\rm (\ref{lpp1})} up to order \( N\) and \( \mbox{\mathversion{bold}$h$}_{\varepsilon }\rightarrow \mbox{\mathversion{bold}$0$}\) in
\( Y^{N}\) as \( \varepsilon \rightarrow 0\).
\end{pr}
{\it Proof}.
We write the equation in (\ref{lpp1}) as
\begin{eqnarray*}
\mbox{\mathversion{bold}$u$}_{t}=- \alpha \varepsilon \mbox{\mathversion{bold}$u$}_{tx} + P(x,t,\partial_{x})\mbox{\mathversion{bold}$u$} + \mbox{\mathversion{bold}$g$} .
\end{eqnarray*}
Setting \( \mbox{\mathversion{bold}$\phi$} _{1} (x):= \mbox{\mathversion{bold}$u$}_{t}(x,0)\) and taking the trace \( t=0\) of the equation we have
\begin{eqnarray}
\alpha \varepsilon \mbox{\mathversion{bold}$\phi$} _{1} ^{\prime} + \mbox{\mathversion{bold}$\phi$} _{1} 
=P(x,0,\partial_{x})\mbox{\mathversion{bold}$u$}_{0}+\mbox{\mathversion{bold}$f$}(x,0)+\mbox{\mathversion{bold}$h$}_{\varepsilon}(x,0).
\label{cc1}
\end{eqnarray}
Using the notations in (\ref{ccq1}) we have  \( P(x,0,\partial_{x})\mbox{\mathversion{bold}$u$}_{0}+\mbox{\mathversion{bold}$f$}(x,0) = \mbox{\mathversion{bold}$Q$}_{1}(x,0) \). 
A prime denotes a derivative with respect to \( x\). As before, solving the above ordinary differential equation for \( \mbox{\mathversion{bold}$\phi$}_{1}\) under the constraint 
$ \displaystyle \lim_{x\rightarrow \infty}\mbox{\mathversion{bold}$\phi$}_{1}(x)=\mbox{\mathversion{bold}$0$}$ we have
\begin{eqnarray*}
\mbox{\mathversion{bold}$\phi$} _{1}(x) = -\frac{1}{\alpha \varepsilon}\int ^{\infty}_{x}{\rm e}^{-\frac{1}{\alpha \varepsilon}(x-y)}
\big\{  {\mbox{\mathversion{bold}$Q$}}_{1}(y,0)+{\mbox{\mathversion{bold}$h$}}_{\varepsilon}(y,0)\big\} {\rm d}y.
\end{eqnarray*}
We give an ansatz for the form of \( {\mbox{\mathversion{bold}$h$}}_{\varepsilon}\), namely
\begin{eqnarray*}
{\mbox{\mathversion{bold}$h$}}_{\varepsilon}(x,t)=\left( \sum ^{N}_{j=0}\mbox{\mathversion{bold}$C$}_{j,\varepsilon }\frac{t^{j}}{j!}\right){\rm e}^{-x},
\end{eqnarray*}
where \( \mbox{\mathversion{bold}$C$}_{j,\varepsilon }\), \( j=0,1,...,N \), are constant vectors depending on \( \varepsilon \) to be determined later.
From Definition \ref{cclpp1} the first order compatibility condition is
\begin{eqnarray*}
\int^{\infty}_{0}{\rm e}^{\frac{y}{\alpha \varepsilon }}\big\{ {\mbox{\mathversion{bold}$Q$}}^{\prime}_{1}(y,0)+{\mbox{\mathversion{bold}$h$}}^{\prime}_{\varepsilon}(y,0)\big\}{\rm d}y
=\mbox{\mathversion{bold}$0$}.
\end{eqnarray*}
Substituting the ansatz for \( \mbox{\mathversion{bold}$h$}_{\varepsilon }(x,t)\), we have
\begin{eqnarray*}
\mbox{\mathversion{bold}$C$}_{0,\varepsilon}\left( 1-\frac{1}{\alpha \varepsilon }\right) ^{-1} = 
\int ^{\infty}_{0}{\rm e}^{\frac{y}{\alpha \varepsilon}}{\mbox{\mathversion{bold}$Q$}}^{\prime}_{1}(y,0){\rm d}y.
\end{eqnarray*}
Since \( {\mbox{\mathversion{bold}$Q$}}^{\prime}_{1}(0,0)={\mbox{\mathversion{bold}$0$}}\) from the compatibility condition for (\ref{lpd1}), we have by integration by parts
\begin{eqnarray*}
\mbox{\mathversion{bold}$C$}_{0,\varepsilon}= \left( \alpha \varepsilon -1 \right) \int ^{\infty}_{0}
{\rm e}^{\frac{y}{\alpha \varepsilon}}{\mbox{\mathversion{bold}$Q$}}^{\prime \prime}_{1}(y,0){\rm d}y.
\end{eqnarray*}
So if we limit ourselves to \( 0< \varepsilon <\min \{1, 1/|\alpha |\} \), from 
\begin{eqnarray*}
{\rm e}^{\frac{y}{\alpha \varepsilon }}|{\mbox{\mathversion{bold}$Q$}}^{\prime \prime}_{1}(y,0)|\leq {\rm e}^{-y}|{\mbox{\mathversion{bold}$Q$}}^{\prime \prime}_{1}(y,0)|,
\end{eqnarray*}
and for \( y>0\)
\begin{eqnarray*}
{\rm e}^{\frac{y}{\alpha \varepsilon }}|{\mbox{\mathversion{bold}$Q$}}^{\prime \prime}_{1}(y,0)| \rightarrow 0 \ {\rm as} \ \varepsilon \rightarrow 0,
\end{eqnarray*}
we see that \( \mbox{\mathversion{bold}$C$}_{0,\varepsilon }\rightarrow \mbox{\mathversion{bold}$0$}\) as \( \varepsilon \rightarrow 0\).
We will show by induction that \( \mbox{\mathversion{bold}$C$}_{j,\varepsilon }\) can be chosen so that 
\( \mbox{\mathversion{bold}$C$}_{j,\varepsilon }\rightarrow \mbox{\mathversion{bold}$0$}\) for \( 1\leq j\leq N\) and \( \mbox{\mathversion{bold}$g$}= 
\mbox{\mathversion{bold}$f$} + \mbox{\mathversion{bold}$h$}_{\varepsilon }\) with \( \mbox{\mathversion{bold}$u$}_{0}\) and \( \mbox{\mathversion{bold}$w$}\)
 satisfies the compatibility conditions for (\ref{lpp1}) up to order \( N\).
Suppose that the above statement holds for \( 0\leq j\leq n-2\) for some \( n\) with \( 2\leq n\leq N\).

We define \( \mbox{\mathversion{bold}$P$}_{n}(x,0)\) and \( \mbox{\mathversion{bold}$\phi$}_{n}(x)\) as in subsection 3.2 and we have
\begin{eqnarray}
\mbox{\mathversion{bold}$\phi$}_{n}(x)=-\frac{1}{\alpha \varepsilon }\int^{\infty}_{x}{\rm e}^{-\frac{1}{\alpha \varepsilon}(x-y)}
\mbox{\mathversion{bold}$P$}_{n}(y,0){\rm d}y,
\label{phin}
\end{eqnarray}
and the \( n\)-th order compatibility condition for (\ref{lpp1}) is
\begin{eqnarray*}
\int^{\infty}_{0}{\rm e}^{\frac{y}{\alpha \varepsilon}}
\mbox{\mathversion{bold}$P$}_{n}^{\prime}(y,0){\rm d}y=\mbox{\mathversion{bold}$0$}.
\end{eqnarray*}
We rewrite this condition as
\begin{eqnarray}
-\mbox{\mathversion{bold}$P$}_{n}^{\prime}(0,0) + \int^{\infty}_{0}{\rm e}^{\frac{y}{\alpha \varepsilon}}
\mbox{\mathversion{bold}$P$}_{n}^{\prime \prime}(y,0){\rm d}y = \mbox{\mathversion{bold}$0$}
\label{cccn}
\end{eqnarray}
by integration by parts. 
We recall that \( \mbox{\mathversion{bold}$P$}_{n}(x,0)\) was successively defined by
\begin{eqnarray*}
\mbox{\mathversion{bold}$P$}_{n}(\cdot,0) = \alpha \mbox{\mathversion{bold}$\phi$}_{n-1}^{\prime \prime \prime}+ \sum ^{n-1}_{j=0}
\left(
\begin{array}{c}
n-1\\
j
\end{array}\right)
{\rm B}_{j}\mbox{\mathversion{bold}$\phi$}_{n-1-j} + \partial ^{n-1}_{t}\mbox{\mathversion{bold}$g$}(\cdot, 0),
\end{eqnarray*}
with \( \mbox{\mathversion{bold}$P$}_{1}(x,0)= \alpha \mbox{\mathversion{bold}$u$}_{0xxx} + {\rm A}(\mbox{\mathversion{bold}$w$}(x ,0),\partial _{x})\mbox{\mathversion{bold}$u$}_{0}
+ \mbox{\mathversion{bold}$g$}(x,0)\).
Substituting (\ref{phin}) with \( n=j\) for \( \mbox{\mathversion{bold}$\phi$}_{j}\) and using integration by parts, we have
\begin{eqnarray*}
\begin{aligned}
\mbox{\mathversion{bold}$P$}_{n}(\cdot ,0)= \alpha \mbox{\mathversion{bold}$P$}_{n-1}^{\prime \prime \prime}
+ \sum^{n-1}_{j=0}\left(
\begin{array}{c}
n-1\\
j
\end{array}\right)
{\rm B}_{j}
&\mbox{\mathversion{bold}$P$}_{n-1-j} + \partial ^{n-1}_{t}\mbox{\mathversion{bold}$g$}(\cdot ,0)\\
& - \alpha \varepsilon 
\bigg\{\alpha \mbox{\mathversion{bold}$\phi$}_{n-1}^{\prime \prime \prime \prime}
+ \sum ^{n-1}_{j=0}
\left(
\begin{array}{c}
n-1\\
j
\end{array}\right)
{\rm B}_{j}\mbox{\mathversion{bold}$\phi$}_{n-1-j}^{\prime}\bigg\}.
\end{aligned}
\end{eqnarray*}
Also recall that 
\begin{eqnarray*}
\mbox{\mathversion{bold}$Q$}_{n}= \alpha \partial^{3}_{x}\mbox{\mathversion{bold}$Q$}_{n-1}
+ \sum^{n-1}_{j=0}\left(
\begin{array}{c}
n-1\\
j
\end{array}\right)
{\rm B}_{j}
\mbox{\mathversion{bold}$Q$}_{n-1-j} + \partial ^{n-1}_{t}\mbox{\mathversion{bold}$f$},
\end{eqnarray*}
with \( \mbox{\mathversion{bold}$Q$}_{1}(x,0) = \alpha \mbox{\mathversion{bold}$u$}_{0xxx} + {\rm A}(\mbox{\mathversion{bold}$w$}(x ,0),\partial _{x})\mbox{\mathversion{bold}$u$}_{0}
+ \mbox{\mathversion{bold}$f$}(x,0)\).
Thus, setting \( \mbox{\mathversion{bold}$R$}_{n} := \mbox{\mathversion{bold}$P$}_{n}-\mbox{\mathversion{bold}$Q$}_{n}\), we have
\begin{eqnarray*}
\begin{aligned}
\mbox{\mathversion{bold}$R$}_{n}(x,0)=\alpha \mbox{\mathversion{bold}$R$}_{n-1}^{\prime \prime \prime}
+ \sum^{n-1}_{j=0}\left(
\begin{array}{c}
n-1\\
j
\end{array}\right)
{\rm B}_{j}
&\mbox{\mathversion{bold}$R$}_{n-1-j} 
+ \partial ^{n-1}_{t}\mbox{\mathversion{bold}$h$}_{\varepsilon }(\cdot ,0) \\
&- \alpha \varepsilon \bigg\{ \alpha \mbox{\mathversion{bold}$\phi$}_{n-1}^{\prime \prime \prime \prime}
+ \sum ^{n-1}_{j=0}
\left(
\begin{array}{c}
n-1\\
j
\end{array}\right)
{\rm B}_{j}\mbox{\mathversion{bold}$\phi$}_{n-1-j}^{\prime}\bigg\},
\end{aligned}
\end{eqnarray*}
with \( \mbox{\mathversion{bold}$R$}_{1}(x,0)= \mbox{\mathversion{bold}$h$}_{\varepsilon }(x,0)\).
We prove by induction that \( \mbox{\mathversion{bold}$R$}_{n}(x,0) = \partial ^{n-1}_{t}\mbox{\mathversion{bold}$h$}_{\varepsilon }(x,0) + o(1)\), 
where \( o(1)\) are terms that tend to zero as \( \varepsilon \rightarrow 0\). The case \( n=1\) is obvious from the definition of 
\( \mbox{\mathversion{bold}$R$}_{1}(x,0)\). Suppose that it holds for \( \mbox{\mathversion{bold}$R$}_{k}(x,0)\) for \( 1\leq k\leq n-1\). 
From the above expression for \( \mbox{\mathversion{bold}$R$}_{n}(x,0)\), the assumption of induction on 
\( \mbox{\mathversion{bold}$R$}_{n}\), and the assumption of induction that 
\( \mbox{\mathversion{bold}$C$}_{j,\varepsilon }\rightarrow \mbox{\mathversion{bold}$0$}\) for \( 0\leq j\leq n-2\), we see that
\begin{eqnarray*}
\mbox{\mathversion{bold}$R$}_{n}(x,0) = \partial ^{n-1}_{t}\mbox{\mathversion{bold}$h$}_{\varepsilon } + o(1)
- \alpha \varepsilon \bigg\{\alpha \mbox{\mathversion{bold}$\phi$}_{n-1}^{\prime \prime \prime \prime}
+ \sum ^{n-1}_{j=0}
\left(
\begin{array}{c}
n-1\\
j
\end{array}\right)
{\rm B}_{j}\mbox{\mathversion{bold}$\phi$}_{n-1-j}^{\prime}\bigg\}.
\end{eqnarray*}
Again, from (\ref{phin}) and Lebesgue's dominated convergence theorem, we see that the last two terms are \( o(1)\), which proves 
\( \mbox{\mathversion{bold}$R$}_{n}(x,0) = \mbox{\mathversion{bold}$P$}_{n}(x,0)-\mbox{\mathversion{bold}$Q$}_{n}(x,0)
	=\partial ^{n-1}_{t}\mbox{\mathversion{bold}$h$}_{\varepsilon }(x,0) + o(1)\). Here, we have used the fact that 
\( \mbox{\mathversion{bold}$P$}_{k}(x,0)\) for \( 1\leq k\leq n-1\) are uniformly bounded with respect to \( \varepsilon \).
We note that from the expressions of \( \mbox{\mathversion{bold}$R$}_{n}(x,0)\) and \( \mbox{\mathversion{bold}$h$}_{\varepsilon }\), 
the terms in \( o(1)\) are composed of terms such that their \( x\) derivative are 
also \( o(1)\). Substituting for \( \mbox{\mathversion{bold}$P$}_{n}(x,0)\) and the ansatz for \( \mbox{\mathversion{bold}$h$}_{\varepsilon }\) in (\ref{cccn}) yields, 
\begin{eqnarray*}
\begin{aligned}
\mbox{\mathversion{bold}$C$}_{n-1,\varepsilon } &= \mbox{\mathversion{bold}$Q$}_{n}^{\prime}(0,0)  
+ \int ^{\infty}_{0}{\rm e}^{\frac{y}{\alpha \varepsilon }}\mbox{\mathversion{bold}$Q$}_{n}^{\prime \prime}(y,0){\rm d}y+  o(1)\\
&=\int ^{\infty}_{0}{\rm e}^{\frac{y}{\alpha \varepsilon }}\mbox{\mathversion{bold}$Q$}_{n}^{\prime \prime}(y,0){\rm d}y+  o(1),
\end{aligned}
\end{eqnarray*}
where we have used the assumption of induction that \( \mbox{\mathversion{bold}$u$}_{0}\), \( \mbox{\mathversion{bold}$f$}\), and 
\( \mbox{\mathversion{bold}$w$}\) satisfy the \( n\)-th order compatibility condition for (\ref{lpd1}), i.e. 
\( \mbox{\mathversion{bold}$Q$}_{n}^{\prime}(0,0)=\mbox{\mathversion{bold}$0$}\).
By using the above expression to define \( \mbox{\mathversion{bold}$C$}_{n-1,\varepsilon }\), we see that 
\( \mbox{\mathversion{bold}$C$}_{n-1,\varepsilon }\rightarrow \mbox{\mathversion{bold}$0$}\) as \( \varepsilon \rightarrow 0\) and 
\( \mbox{\mathversion{bold}$u$}_{0}\), \( \mbox{\mathversion{bold}$g$}\), and \( \mbox{\mathversion{bold}$w$}\) satisfy the compatibility condition 
for (\ref{lpp1}) up to order \( n\). 
Furthermore, from the explicit form we see that \( \mbox{\mathversion{bold}$h$}_{\varepsilon }\rightarrow \mbox{\mathversion{bold}$0$}\) in \( Y^{N}\). 
This finishes the proof of the proposition.
\hfill \( \Box\)
\bigskip

The corrections to the data associated with (\ref{lp1}) can be treated the same way.
 
%
%
%
%
\section{Construction and Estimate of Solution for the Regularized System}
\setcounter{equation}{0}
We construct the solution to (\ref{lp1}) as a sum of two functions \( \mbox{\mathversion{bold}$u$}_{1}\) and 
\( \mbox{\mathversion{bold}$u$}_{2}\) which are defined as the solutions of the following systems. \( \mbox{\mathversion{bold}$u$}_{1}\) is defined as the solution
to the initial value problem
\begin{eqnarray}
\left\{
\begin{array}{ll}
\mbox{\mathversion{bold}$u$}_{1t}= \alpha \left( \mbox{\mathversion{bold}$u$}_{1xx} - \varepsilon \mbox{\mathversion{bold}$u$}_{1t}\right)_{x} + \mbox{\mathversion{bold}$G$}
, & x \in {\mathbf R}, t>0, \\
\mbox{\mathversion{bold}$u$}_{1}(x,0) = \mbox{\mathversion{bold}$U$}_{0}, & x\in {\mathbf R}, \\
\end{array}\right.
\label{one}
\end{eqnarray}
and \( \mbox{\mathversion{bold}$u$}_{2}\) is defined as the solution to the initial-boundary value problem
\begin{eqnarray}
\left\{
\begin{array}{ll}
\mbox{\mathversion{bold}$u$}_{2t}= \alpha \left( \mbox{\mathversion{bold}$u$}_{2xx}-\varepsilon \mbox{\mathversion{bold}$u$}_{2t}\right)_{x}, & x>0, t>0, \\
\mbox{\mathversion{bold}$u$}_{2}(x,0)= \mbox{\mathversion{bold}$0$}, & x>0, \\
\mbox{\mathversion{bold}$u$}_{2x}(0,t)=-\mbox{\mathversion{bold}$u$}_{1x}(0,t)=:\mbox{\mathversion{bold}$\Phi$}(t),& t>0.
\end{array}\right.
\label{two1}
\end{eqnarray}
Here, \( \mbox{\mathversion{bold}$G$}\) and \( \mbox{\mathversion{bold}$U$}_{0}\) are smooth extensions of \( \mbox{\mathversion{bold}$g$} \) and \( \mbox{\mathversion{bold}$u$}_{0}\)
to \( x<0\), respectively.

\subsection{Construction and Estimate of \( \mbox{\mathversion{bold}$u$}_{1}\) }

First we solve (\ref{one}). By taking the Fourier transform with respect to \( x\), we obtain the ordinary differential equation
\begin{eqnarray*}
\left\{
\begin{array}{ll}
\hat{\mbox{\mathversion{bold}$u$}}_{1t}=
\frac{1}{1+{\rm i}\alpha \varepsilon \xi}\left( -{\rm i}\alpha \xi ^{3} \hat{\mbox{\mathversion{bold}$u$}}_{1} + \hat{\mbox{\mathversion{bold}$G$}}\right) ,\\
\hat{\mbox{\mathversion{bold}$u$}}_{1} (\xi ,0)=\hat{\mbox{\mathversion{bold}$U$}}_{0},
\end{array}\right.
\end{eqnarray*}
where \( \hat{\mbox{\mathversion{bold}$u$}}_{1}\) is the Fourier transform defined by
\begin{eqnarray*}
\hat{\mbox{\mathversion{bold}$u$}}_{1}(\xi ,t) = \frac{1}{\sqrt{2\pi} }\int ^{\infty}_{-\infty}{\rm e}^{-{\rm i}x\xi }\mbox{\mathversion{bold}$u$}_{1}(x,t){\rm d}x.
\end{eqnarray*}
The solution can be explicitly constructed as
\begin{eqnarray*}
\hat{{\mbox{\mathversion{bold}$u$}}}_{1}(\xi ,t) = {\rm e}^{c(\xi )t}\hat{{\mbox{\mathversion{bold}$U$}}}_{0}+\int^{t}_{0}{\rm e}^{c(\xi )(t-\tau )}
\frac{1}{1+{\rm i}\alpha \varepsilon \xi }\hat{{\mbox{\mathversion{bold}$G$}}}(\xi ,\tau ){\rm d}\tau ,
\end{eqnarray*}
where \( c(\xi )\) is given by
\begin{eqnarray*}
c(\xi )=\frac{-\alpha ^{2}\varepsilon \xi ^{4}-{\rm i}\alpha \xi ^{3}}{1+\alpha ^{2}\varepsilon ^{2}\xi ^{2}}.
\end{eqnarray*}
Now we derive an estimate for \( {\mbox{\mathversion{bold}$u$}}_{1}\). The estimate derived here will be of parabolic nature, and will not be 
uniform in \( \varepsilon \). 

\begin{eqnarray*}
\begin{aligned}
\frac{1}{2}\frac{{\rm d}}{{\rm d}t}\bigg\{ \| {\mbox{\mathversion{bold}$u$}}_{1}\|_{L^{2}(\mathbf{R})} ^{2} + \alpha ^{2}\varepsilon ^{2}\| {\mbox{\mathversion{bold}$u$}}_{1x}\|
_{L^{2}(\mathbf{R})} ^{2} \bigg\}
&=\langle{\mbox{\mathversion{bold}$u$}}_{1},{\mbox{\mathversion{bold}$u$}}_{1t}\rangle + \alpha ^{2}\varepsilon ^{2}\langle{\mbox{\mathversion{bold}$u$}}_{1x},{\mbox{\mathversion{bold}$u$}}_{1xt}\rangle \\
&\hspace*{-3cm}\leq \frac{1}{2}\bigg( \|{\mbox{\mathversion{bold}$u$}}_{1}\|_{L^{2}(\mathbf{R})}^{2} + \alpha ^{2}\varepsilon ^{2}\| {\mbox{\mathversion{bold}$u$}}_{1x}\|
_{L^{2}(\mathbf{R})} ^{2}\bigg) 
-\alpha ^{2}\varepsilon \| {\mbox{\mathversion{bold}$u$}}_{1xx}\|_{L^{2}(\mathbf{R})} ^{2} + \| {\mbox{\mathversion{bold}$G$}}\|_{L^{2}(\mathbf{R})} ^{2}.
\end{aligned}
\end{eqnarray*}
Similarly, for an integer \( l\) we have
\begin{eqnarray*}
\begin{aligned}
\frac{1}{2}\frac{{\rm d}}{{\rm d}t}\bigg\{ \| \partial ^{l}_{x}{\mbox{\mathversion{bold}$u$}}_{1}\|_{L^{2}(\mathbf{R})} ^{2}
&+ \alpha ^{2}\varepsilon ^{2}\| \partial ^{l+1}_{x}{\mbox{\mathversion{bold}$u$}}_{1}\|_{L^{2}(\mathbf{R})} ^{2}\bigg\} \\
&\hspace*{-2cm}\leq
\frac{1}{2}\bigg(  \| \partial ^{l}_{x}{\mbox{\mathversion{bold}$u$}}_{1}\|_{L^{2}(\mathbf{R})} ^{2}
+ \alpha ^{2}\varepsilon ^{2}\| \partial ^{l+1}_{x}{\mbox{\mathversion{bold}$u$}}_{1}\|_{L^{2}(\mathbf{R})} ^{2}\bigg) 
-\alpha ^{2} \varepsilon \| \partial ^{l+2}_{x}{\mbox{\mathversion{bold}$u$}}_{1}\|_{L^{2}(\mathbf{R})} ^{2}+ \| \partial ^{l}{\mbox{\mathversion{bold}$G$}}\|_{L^{2}(\mathbf{R})} ^{2}.
\end{aligned}
\end{eqnarray*}
We also give estimates for the \( t\) derivatives of \( {\mbox{\mathversion{bold}$u$}}_{1}\) because it will come in use later. 

The estimates are the same as above and 
we obtain for integers \( l\) and \( m\),
\begin{eqnarray*}
\begin{aligned}
\frac{1}{2}\frac{{\rm d}}{{\rm d}t}\bigg\{ \| \partial ^{l}_{x}\partial ^{m}_{t}{\mbox{\mathversion{bold}$u$}}_{1}\|_{L^{2}(\mathbf{R})}^{2}
&+ \alpha ^{2}\varepsilon ^{2}\| \partial ^{l+1}_{x}\partial ^{m}_{t}{\mbox{\mathversion{bold}$u$}}_{1}\|_{L^{2}(\mathbf{R})} ^{2}\bigg\} \\
&\hspace*{-3.3cm}\leq
\frac{1}{2}\bigg( \| \partial ^{l}_{x}\partial ^{m}_{t}{\mbox{\mathversion{bold}$u$}}_{1}\|_{L^{2}(\mathbf{R})}^{2}
+ \alpha ^{2}\varepsilon ^{2}\| \partial ^{l+1}_{x}\partial ^{m}_{t}{\mbox{\mathversion{bold}$u$}}_{1}\|_{L^{2}(\mathbf{R})} ^{2} \bigg) 
- \alpha ^{2} \varepsilon \| \partial ^{l+2}_{x}\partial ^{m}_{t}{\mbox{\mathversion{bold}$u$}}\|_{L^{2}(\mathbf{R})}^{2}
+ \| \partial ^{l}_{x}\partial ^{m}_{t}{\mbox{\mathversion{bold}$G$}}\|_{L^{2}(\mathbf{R})}^{2}.
\end{aligned}
\end{eqnarray*}
To finish this estimate, we must estimate the \( t\) derivatives of \( {\mbox{\mathversion{bold}$u$}}_{1}\) at \( t=0\) in terms of 
\( {\mbox{\mathversion{bold}$U$}}_{0}\) and \( {\mbox{\mathversion{bold}$G$}}\).
Set \( {\mbox{\mathversion{bold}$\phi $}} _{1n}(x) := \partial ^{n}_{t}{\mbox{\mathversion{bold}$u$}}_{1}(x,0) \). As before, by taking the trace \( t=0\) of the equation,
solving for \( \mbox{\mathversion{bold}$\phi $}_{11}\) under the constraint 
$ \displaystyle\lim _{x \to \infty} {\mbox{\mathversion{bold}$\phi $}} _{11}(x) ={\mbox{\mathversion{bold}$0 $}} $
\begin{eqnarray*}
{\mbox{\mathversion{bold}$\phi $}} _{11}(x) = -\frac{1}{\alpha \varepsilon}\int ^{\infty}_{x}{\rm e}^{-\frac{1}{\alpha \varepsilon}(x-y)}
\big\{ \alpha {\mbox{\mathversion{bold}$U$}}_{0}^{\prime \prime \prime}(y)+{\mbox{\mathversion{bold}$G$}}(y,0)\big\} {\rm d}y.
\end{eqnarray*}
Through direct calculation, we see that 
\begin{eqnarray*}
\partial ^{k}_{x}{\mbox{\mathversion{bold}$\phi $}} _{11}(x) = -\frac{1}{\alpha \varepsilon}\int ^{\infty}_{x}{\rm e}^{-\frac{1}{\alpha \varepsilon}(x-y)}
\big\{ \alpha \partial ^{k+3}_{y}{\mbox{\mathversion{bold}$U$}}_{0}(y)+\partial^{k}_{y}{\mbox{\mathversion{bold}$G$}}(y,0)\big\} {\rm d}y.
\end{eqnarray*}
Also from direct calculation, we obtain
\begin{eqnarray}
\| \partial ^{k}_{x}{\mbox{\mathversion{bold}$\phi $}}_{11} \|_{L^{2}(\mathbf{R})} \leq C\| \partial ^{k+2}_{x}{\mbox{\mathversion{bold}$U$}}_{0} 
\|_{L^{2}(\mathbf{R})} + \| \partial ^{k}_{x}{\mbox{\mathversion{bold}$G$}}(\cdot ,0)\|_{L^{2}(\mathbf{R})} .
\label{est}
\end{eqnarray}
Here, we have used integration by parts, that is
\begin{eqnarray*}
\begin{aligned}
&-\frac{1}{\varepsilon}\int ^{\infty}_{x}{\rm e}^{-\frac{1}{\alpha \varepsilon}(x-y)}\partial ^{k+3}_{y}{\mbox{\mathversion{bold}$U$}}_{0}(y) {\rm d}y\\
&\hspace*{3cm}=-\frac{1}{\varepsilon }\left[ {\rm e}^{-\frac{1}{\alpha \varepsilon}(x-y)}\partial ^{k+2}_{y}{\mbox{\mathversion{bold}$U$}}_{0}(y)\right] ^{\infty}_{y=x}
+\frac{1}{\alpha \varepsilon ^{2}}\int ^{\infty}_{x}{\rm e}^{-\frac{1}{\alpha \varepsilon}(x-y)}\partial ^{k+2}_{y}{\mbox{\mathversion{bold}$U$}}_{0}(y){\rm d}y \\
&\hspace*{3cm}= \frac{1}{\varepsilon }\partial ^{k+2}_{x}{\mbox{\mathversion{bold}$U$}}_{0}(x)
+ \int ^{\infty}_{x}{\rm e}^{-\frac{1}{\alpha \varepsilon}(x-y)}\partial ^{k+2}_{y}{\mbox{\mathversion{bold}$U$}}_{0}(y){\rm d}y.
\end{aligned}
\end{eqnarray*}
As shown from the above calculation, the constant \( C\) in (\ref{est}) is not uniform in \( \varepsilon \).

Taking the \( t\) derivative of the equation \( n-1\) times, we obtain
\begin{eqnarray*}
{\mbox{\mathversion{bold}$\phi $}} _{1n}^{\prime}= \frac{1}{\alpha \varepsilon}\big\{ -{\mbox{\mathversion{bold}$\phi $}} _{1n} 
+ \alpha \phi _{1(n-1)}^{\prime \prime \prime}+ \partial^{n-1}_{t}{\mbox{\mathversion{bold}$G$}}(\cdot ,0) \big\}.
\end{eqnarray*}
As before, we obtain the following expression and estimate
\begin{eqnarray*}
\partial^{k}_{x}{\mbox{\mathversion{bold}$\phi $}} _{1n}(x) = -\frac{1}{\alpha \varepsilon}\int ^{\infty}_{x}{\rm e}^{-\frac{1}{\alpha \varepsilon}(x-y)}
\big\{ \alpha \partial ^{k+3}_{y}{\mbox{\mathversion{bold}$\phi $}} _{1(n-1)}(y)+\partial^{n-1}_{t}\partial^{k}_{y}{\mbox{\mathversion{bold}$G$}}(y,0)\big\} {\rm d}y,
\end{eqnarray*}
\begin{eqnarray*}
\| \partial^{k}_{x}\mbox{\mathversion{bold}$\phi $} _{1n}\|_{L^{2}(\mathbf{R})} \leq C\bigg( \| {\mbox{\mathversion{bold}$U$}}_{0}\|_{H^{k+2n}(\mathbf{R})}
+ \sum ^{n-1}_{j=0} \| \partial ^{n-1-j}_{t}{\mbox{\mathversion{bold}$G$}}(\cdot ,0)\|_ {H^{k+2j}(\mathbf{R})} \bigg),
\end{eqnarray*}
where \( C\) is a positive constant depending on \( \varepsilon \). Combining these estimates yields
\begin{eqnarray*}
\begin{aligned}
\sup_{0\leq t \leq T} \bigg\{ \| \partial ^{m}_{t}\partial^{l}_{x}{\mbox{\mathversion{bold}$u$}}_{1}(\cdot, t)\|_{L^{2}(\mathbf{R})} ^{2} 
&+ \alpha^{2}\varepsilon^{2}\| \partial ^{m}_{t}\partial^{l+1}_{x}{\mbox{\mathversion{bold}$u$}}_{1}(\cdot, t)\|_{L^{2}(\mathbf{R})}^{2} \bigg\}
+ \int^{T}_{0}\| \partial^{m}_{t}\partial^{l+2}_{x}{\mbox{\mathversion{bold}$u$}}_{1}(\cdot ,t )\|_{L^{2}(\mathbf{R})}^{2}{\rm d}t \\
&\hspace*{-4cm}\leq C{\rm e}^{T}\bigg( \| {\mbox{\mathversion{bold}$U$}}_{0}\|^{2}_{H^{l+2m}(\mathbf{R})}
+ \sum ^{n-1}_{j=0} \| \partial ^{n-1-j}_{t}{\mbox{\mathversion{bold}$G$}}(\cdot ,0)\|^{2}_ {H^{l+2j}(\mathbf{R})}\bigg) 
+ C\int ^{T}_{0}{\rm e}^{T-t}\| \partial ^{m}_{t}\partial ^{l}_{x}{\mbox{\mathversion{bold}$G$}}(\cdot, t)\|_{L^{2}(\mathbf{R})}^{2}{\rm d}t.
\end{aligned}
\end{eqnarray*}
From the boundedness of the extension, we have the following estimate on the half-line.
\begin{eqnarray*}
\begin{aligned}
\sup_{0\leq t \leq T} \bigg\{ \| \partial ^{m}_{t}\partial^{l}_{x}{\mbox{\mathversion{bold}$u$}}_{1}(\cdot, t)\| ^{2} 
&+ \alpha^{2}\varepsilon^{2}\| \partial ^{m}_{t}\partial^{l+1}_{x}{\mbox{\mathversion{bold}$u$}}_{1}(\cdot, t)\| ^{2} \bigg\}
+ \int^{T}_{0}\| \partial^{m}_{t}\partial^{l+2}_{x}{\mbox{\mathversion{bold}$u$}}_{1}(\cdot ,t )\| ^{2}{\rm d}t \\
&\hspace*{-2.5cm}\leq C{\rm e}^{T}\bigg( \| {\mbox{\mathversion{bold}$u$}}_{0}\|^{2}_{l+2m}
+ \sum ^{n-1}_{j=0} \| \partial ^{n-1-j}_{t}{\mbox{\mathversion{bold}$g$}}(\cdot ,0)\|^{2}_ {l+2j}\bigg) 
+ C\int ^{T}_{0}{\rm e}^{T-t}\| \partial ^{m}_{t}\partial ^{l}_{x}{\mbox{\mathversion{bold}$g$}}(\cdot, t)\| ^{2}{\rm d}t.
\end{aligned}
\end{eqnarray*}

\subsection{Construction and Estimate of \( \mbox{\mathversion{bold}$u$}_{2}\) }

Now we solve (\ref{two1}). First we take a look at the compatibility condition. We derive the compatibility conditions for (\ref{two1}) and check that it is
satisfied at the same time.
Suppose that the initial data and the forcing term satisfy the compatibility conditions for (\ref{lp1}) up to some finite order. 
The zero-th order compatibility condition for (\ref{two1}) is \( {\mbox{\mathversion{bold}$u$}}_{1x}(0,0)={\mbox{\mathversion{bold}$0$}}\). From 
the way we defined \( {\mbox{\mathversion{bold}$u$}}_{1}\) and the compatibility condition for (\ref{lp1}), we have
\begin{eqnarray*}
-{\mbox{\mathversion{bold}$u$}}_{1x}(0,0) = -{\mbox{\mathversion{bold}$u$}}_{0x}(0,0) ={\mbox{\mathversion{bold}$0$}},
\end{eqnarray*}
so that the zero-th order compatibility condition for (\ref{two1}) is satisfied. Now we look at the first order compatibility condition. Taking the 
\( t\) derivative of the boundary condition, we have \( {\mbox{\mathversion{bold}$u$}}_{2tx}(0,0)=-{\mbox{\mathversion{bold}$u$}}_{1tx}(0,0) \).
Taking the trace \( t=0\) of the equation in (\ref{two1}) and setting \( {\mbox{\mathversion{bold}$\phi $}}_{21}(x):= {\mbox{\mathversion{bold}$u$}}_{2t}(x,0)\) yield
\begin{eqnarray*}
{\mbox{\mathversion{bold}$\phi $}}_{21}^{\prime} = -\frac{1}{\alpha \varepsilon }{\mbox{\mathversion{bold}$\phi $}}_{21}.
\end{eqnarray*}
Solving for \( {\mbox{\mathversion{bold}$\phi $}}_{21}\), we have
\begin{eqnarray*}
{\mbox{\mathversion{bold}$\phi $}}_{21}(x)= {\mbox{\mathversion{bold}$\phi $}}_{21}(0){\rm e}^{-\frac{1}{\alpha \varepsilon }x}.
\end{eqnarray*}
For \( {\mbox{\mathversion{bold}$\phi $}}_{21}\) to be integrable, \( {\mbox{\mathversion{bold}$\phi $}}_{21}(0)\) must be a zero vector, thus making 
\( {\mbox{\mathversion{bold}$\phi $}}_{21}(x) = {\mbox{\mathversion{bold}$0$}}\) for any \( x>0\), from which we can deduce that the first order
compatibility condition for (\ref{two1}) is \( {\mbox{\mathversion{bold}$u$}}_{1tx}(0,0)=
{\mbox{\mathversion{bold}$0$}}\).
Taking the trace \( t=0\) of the equation in (\ref{one}) and setting \( {\mbox{\mathversion{bold}$\phi $}}_{11}(x):= {\mbox{\mathversion{bold}$u$}}_{1t}(x,0) \) we have
\begin{eqnarray*}
{\mbox{\mathversion{bold}$\phi $}}_{11}^{\prime} = -\frac{1}{\alpha \varepsilon }{\mbox{\mathversion{bold}$\phi $}}_{11}
+\frac{1}{\alpha \varepsilon }\big\{ \alpha {\mbox{\mathversion{bold}$U$}}^{\prime \prime \prime}_{0} + {\mbox{\mathversion{bold}$G$}} \big\}.
\end{eqnarray*}
As before, solving for \( {\mbox{\mathversion{bold}$\phi $}}_{11}\) and imposing the integrability of \( {\mbox{\mathversion{bold}$\phi $}}_{11}\) we arrive at
\begin{eqnarray*}
{\mbox{\mathversion{bold}$\phi $}}_{11}(x)= -\frac{1}{\alpha \varepsilon }\int ^{\infty}_{x}{\rm e}^{-\frac{1}{\alpha \varepsilon }(x-y)}
\big\{ \alpha {\mbox{\mathversion{bold}$U$}}^{\prime \prime \prime}_{0}(y) + {\mbox{\mathversion{bold}$G$}}(y,0) \big\} {\rm d}y.
\end{eqnarray*}
If \( x\) is restricted to \( x>0\), \( {\mbox{\mathversion{bold}$U$}}_{0}\) and \( {\mbox{\mathversion{bold}$G$}}\) can be replaced with 
\( {\mbox{\mathversion{bold}$u$}}_{0}\) and \( {\mbox{\mathversion{bold}$g$}}\), respectively, because they are extensions of the respective functions.
Taking the trace \( t=0\) of the equation in (\ref{lp1}), setting \( {\mbox{\mathversion{bold}$\phi $}}_{1}(x):= {\mbox{\mathversion{bold}$u$}}_{t}(x,0)\), and
solving for \( {\mbox{\mathversion{bold}$\phi $}}_{1}\) we have
\begin{eqnarray}
\begin{aligned}
{\mbox{\mathversion{bold}$\phi $}}_{1}(x)&=-\frac{1}{\alpha \varepsilon }\int ^{\infty}_{x}{\rm e}^{-\frac{1}{\alpha \varepsilon }(x-y)}
\big\{ \alpha {\mbox{\mathversion{bold}$u$}}^{\prime \prime \prime}_{0}(y) + {\mbox{\mathversion{bold}$g$}}(y,0) \big\} {\rm d}y \\
&= {\mbox{\mathversion{bold}$\phi $}}_{11}(x).
\label{phi1}
\end{aligned}
\end{eqnarray}
Taking the \( t\) derivative of the boundary condition in (\ref{lp1}) and taking the trace \( x=0\) and \( t=0\) we see that 
\( {\mbox{\mathversion{bold}$\phi $}}^{\prime}_{1}(0)={\mbox{\mathversion{bold}$u$}}_{tx}(0,0)={\mbox{\mathversion{bold}$0$}}\), which gives
\begin{eqnarray*}
{\mbox{\mathversion{bold}$u$}}_{1tx}(0,0) &= {\mbox{\mathversion{bold}$\phi $}}^{\prime}_{11}(0) 
= {\mbox{\mathversion{bold}$\phi $}}^{\prime}_{1}(0)
= {\mbox{\mathversion{bold}$0$}},
\end{eqnarray*}
where we have used (\ref{phi1}). The above shows that the first order compatibility condition for (\ref{two1}) is satisfied. 

In the same manner, 
we set \( {\mbox{\mathversion{bold}$\phi $}}_{1n}(x):= \partial ^{n}_{t}{\mbox{\mathversion{bold}$u$}}_{1}(x,0)\),
\( {\mbox{\mathversion{bold}$\phi $}}_{2n}(x):= \partial ^{n}_{t}{\mbox{\mathversion{bold}$u$}}_{2}(x,0)\), and 
\( {\mbox{\mathversion{bold}$\phi $}}_{n}(x):= \partial ^{n}_{t}{\mbox{\mathversion{bold}$u$}}(x,0)\),
where \( \mbox{\mathversion{bold}$\phi $}_{2n}\) and \( \mbox{\mathversion{bold}$\phi $}_{n}\) can be expressed using given data only as in section 3.
We will show that the \( n\)-th order compatibility condition for (\ref{two1}) is satisfied by proving that 
\( {\mbox{\mathversion{bold}$\phi $}}_{1n}={\mbox{\mathversion{bold}$\phi $}}_{n}\) and \({\mbox{\mathversion{bold}$\phi $}}_{2n}=
{\mbox{\mathversion{bold}$0$}}\). We prove this by induction, so suppose that 
\( {\mbox{\mathversion{bold}$\phi $}}_{1k}={\mbox{\mathversion{bold}$\phi $}}_{k}\) and \( {\mbox{\mathversion{bold}$\phi $}}_{2k}={\mbox{\mathversion{bold}$0$}}\) 
for \( k=1,2, \dots n-1\).  
We note that from the compatibility conditions for (\ref{lp1}), \( \mbox{\mathversion{bold}$\phi $}^{\prime}_{k}(0) = \mbox{\mathversion{bold}$0$}\) for \( 0\leq k\leq n\).
These functions satisfy
\begin{eqnarray*}
\begin{aligned}
{\mbox{\mathversion{bold}$\phi $}}_{1n}&= \alpha {\mbox{\mathversion{bold}$\phi $}}_{1(n-1)}^{\prime \prime \prime}-\alpha \varepsilon {\mbox{\mathversion{bold}$\phi $}}_{1n}^{\prime}
+\partial ^{n-1}_{t}{\mbox{\mathversion{bold}$G$}}, \\
{\mbox{\mathversion{bold}$\phi $}}_{2n}&= \alpha {\mbox{\mathversion{bold}$\phi $}}_{2(n-1)}^{\prime \prime \prime} -\alpha \varepsilon {\mbox{\mathversion{bold}$\phi $}}_{2n}^{\prime},\\
{\mbox{\mathversion{bold}$\phi $}}_{n} &= \alpha {\mbox{\mathversion{bold}$\phi $}}_{n-1}^{\prime \prime \prime} -\alpha \varepsilon {\mbox{\mathversion{bold}$\phi $}}_{n}^{\prime}
+\partial ^{n-1}_{t}{\mbox{\mathversion{bold}$g$}}.
\end{aligned}
\end{eqnarray*}
The above is obtained by taking the derivative of the respective equations \( n-1\) times with respect to \( t\) and taking the trace \( t=0\). 
First, we see from \( {\mbox{\mathversion{bold}$\phi $}}_{2(n-1)}={\mbox{\mathversion{bold}$0$}}\) that
\begin{eqnarray*}
{\mbox{\mathversion{bold}$\phi $}}_{2n}^{\prime} = -\frac{1}{\alpha \varepsilon }{\mbox{\mathversion{bold}$\phi $}}_{2n}.
\end{eqnarray*}
As before, from the above equation along with the necessity that \( {\mbox{\mathversion{bold}$\phi $}}_{2n}\) be integrable, we see that \( {\mbox{\mathversion{bold}$\phi $}}_{2n}
={\mbox{\mathversion{bold}$0$}}\). This implies that, through the boundary condition, the \( n\)-th order compatibility condition for (\ref{two1}) is 
\( \partial ^{n}_{t}{\mbox{\mathversion{bold}$u$}}_{1x}(0,0)={\mbox{\mathversion{bold}$0$}}\).
 Solving for \( {\mbox{\mathversion{bold}$\phi $}}_{1n}\) and \( {\mbox{\mathversion{bold}$\phi $}}_{n}\) we have
\begin{eqnarray*}
\begin{aligned}
{\mbox{\mathversion{bold}$\phi $}}_{1n}(x)&=-\frac{1}{\alpha \varepsilon }\int ^{\infty}_{x}{\rm e}^{-\frac{1}{\alpha \varepsilon }(x-y)}\big\{
\alpha {\mbox{\mathversion{bold}$\phi $}}_{1(n-1)}^{\prime \prime \prime}(y)+\partial ^{n-1}_{t}{\mbox{\mathversion{bold}$G$}}(y,0)\big\} {\rm d}y, \\
{\mbox{\mathversion{bold}$\phi $}}_{n}(x)&=-\frac{1}{\alpha \varepsilon }\int ^{\infty}_{x}{\rm e}^{-\frac{1}{\alpha \varepsilon }(x-y)}\big\{
\alpha {\mbox{\mathversion{bold}$\phi $}}_{n-1}^{\prime \prime \prime}(y)+\partial ^{n-1}_{t}{\mbox{\mathversion{bold}$g$}}(y,0)\big\} {\rm d}y.
\end{aligned}
\end{eqnarray*}
Again, from the assumption of induction and the fact that \( {\mbox{\mathversion{bold}$U$}}_{0}\) and \( {\mbox{\mathversion{bold}$G$}}\) are 
extensions of \( {\mbox{\mathversion{bold}$u$}}_{0}\) and \( {\mbox{\mathversion{bold}$g$}}\), respectively, we see that 
\( {\mbox{\mathversion{bold}$\phi $}}_{1n}(x) = {\mbox{\mathversion{bold}$\phi $}}_{n}(x)\). So we have 
\begin{eqnarray*}
\partial ^{n}_{t}{\mbox{\mathversion{bold}$u$}}_{1x}(0,0)=\mbox{\mathversion{bold}$\phi $}^{\prime}_{1n}(0)=\mbox{\mathversion{bold}$\phi $}^{\prime}_{n}(0)
= {\mbox{\mathversion{bold}$0$}},
\end{eqnarray*}
which shows that the \( n\)-th order compatibility condition for (\ref{two1}) is satisfied. 

Now we construct \( {\mbox{\mathversion{bold}$u$}}_{2}\). From the preceding arguments, we saw that 
\( \frac{{\rm d}^{k}\mbox{\mathversion{bold}$\Phi$}}{{\rm d}t^{k}} (0) = \partial ^{k}_{t}\mbox{\mathversion{bold}$u$}_{1x}(0,0) =\mbox{\mathversion{bold}$0$} \) for
\( 0\leq k\leq n\),
so it is natural to construct and estimate \( \mbox{\mathversion{bold}$u$}_{2}\) in Sobolev--Slobodetski\u\i \ spaces.
Taking the Laplace transform of the equation yields
\begin{eqnarray*}
\left\{ 
\begin{array}{ll}
\tau \tilde{{\mbox{\mathversion{bold}$u$}}}_{2}= \alpha \tilde{{\mbox{\mathversion{bold}$u$}}}_{2xxx}-\alpha \varepsilon \tau \tilde{{\mbox{\mathversion{bold}$u$}}}_{2x}, &
x>0, \\
\tilde{{\mbox{\mathversion{bold}$u$}}}_{2x}(0, \tau )=-\tilde{{\mbox{\mathversion{bold}$u$}}}_{1x}(0,\tau ) = \tilde{{\mbox{\mathversion{bold}$\Phi$}}}(\tau ), & \
\end{array}\right.
\end{eqnarray*}
where \( \tau = h+{\rm i}\eta \) with \( h>0\) and \( \eta \in \mathbf{R}\). We show the following properties about the characteristic roots of the above 
ordinary differential equation.
\begin{lm}
For \( h >0 \) and \( \varepsilon >0\), the characteristic equation, \( \lambda ^{3}-\varepsilon \tau \lambda -\frac{\tau}{\alpha } =0 \), has exactly one
root \( \lambda \) satisfying \( \Re \lambda <0\). We will denote this root as \( \mu \). Furthermore, there are positive constants \( \eta _{0}\) and \( C\) such that for 
\( |\eta |\geq \eta _{0}\) the following holds.
\begin{eqnarray*}
\left| \mu - \sqrt{\frac{\varepsilon}{2}}\left( -1+{\rm i}\right)|\eta |^{1/2} \right| \leq C.
\end{eqnarray*}
\label{root}
\end{lm}

\medskip
\noindent
{\it Proof.}
First, we look at the asymptotic behavior of the roots as \( \eta \rightarrow +\infty\). Dividing the characteristic equation by \( \eta ^{3/2}\) and setting 
\( \tilde{\lambda} := \frac{\lambda }{\eta ^{1/2}}\) we have 
\begin{eqnarray}
\tilde{\lambda}^{3}-\frac{\varepsilon h}{\eta}\tilde{\lambda}-{\rm i}\varepsilon \tilde{\lambda} -\frac{h}{\alpha \eta ^{3/2}}-{\rm i}\frac{1}{\alpha \eta ^{1/2}}=0.
\label{asy}
\end{eqnarray}
Taking the limit \( \eta \rightarrow +\infty\), we have
\begin{eqnarray*}
\tilde{\lambda}^{3}-{\rm i}\varepsilon \tilde{\lambda}=0.
\end{eqnarray*}
The solutions are \( \tilde{\lambda}= 0, \ \pm \sqrt{\frac{\varepsilon}{2}}(-1+{\rm i}) \). 
The root \( \sqrt{\frac{\varepsilon}{2}}(-1+{\rm i}) \) corresponds to our desired root of the original characteristic equation.
We must consider the root \( 0\) in more detail. By setting 
\( \tilde{\lambda}=0+ c_{1}\eta ^{-1/2}+{\rm O}(\eta ^{-1})\) and substituting it into (\ref{asy}), we get from the coefficients of terms \( {\rm O}(\eta ^{-1/2}) \),
\begin{eqnarray*}
-{\rm i}\varepsilon c_{1}-{\rm i}\frac{1}{\alpha }=0.
\end{eqnarray*}
This gives \( c_{1}=-\frac{1}{\alpha \varepsilon } >0 \), hence only one root with a negative real part exists for sufficiently large \( \eta \). 
The case \( \eta \rightarrow -\infty \) can be 
treated in the same way. Now we show that for any \( h>0\) and \( \eta \in \mathbf{R}\), there are no pure imaginary roots, which, combined with the continuity of the roots with
respect to the coefficient of the characteristic equation, proves that the number of roots with a negative real part does not change. 

We divide the characteristic equation into its real and imaginary parts. Setting \( \lambda = a+{\rm i}b\) we have
\begin{eqnarray*}
a^{3}-3ab^{2}-\varepsilon ha + \varepsilon \eta b-\frac{h}{\alpha }=0,\\
-b^{3}+3a^{2}b-\varepsilon hb- \varepsilon \eta a-\frac{\eta }{\alpha }=0.
\end{eqnarray*}
Suppose that a pure imaginary root exists, which corresponds to a root with \( a=0\), we then have
\begin{eqnarray*}
\varepsilon \eta b=\frac{h}{\alpha }, \ 
-b^{3}-\varepsilon hb-\frac{\eta }{\alpha }=0.
\end{eqnarray*}
From the first equation we have \( \eta b=\frac{h}{\alpha \varepsilon }\). Multiplying the second equation by \( b\) and substituting for \( \eta b\) yields
\begin{eqnarray*}
-b^{4}-\varepsilon hb^{2}-\frac{h}{\alpha ^{2}\varepsilon }=0.
\end{eqnarray*}
Since we are considering \( h>0\) and \( \varepsilon >0\), the above relation is a contradiction. Thus, no such root exists. 
\hfill \( \Box \)

\bigskip
From Lemma \ref{root}, we see that the Laplace transform of a square integrable solution to (\ref{two1}) can be expressed as
\begin{eqnarray*}
\tilde{{\mbox{\mathversion{bold}$u$}}}_{2}(x,\tau ) = \frac{1}{\mu }\tilde{{\mbox{\mathversion{bold}$\Phi $}}}(\tau ) {\rm e}^{\mu x},
\end{eqnarray*}
where \( \mu \) is the root of the characteristic equation mentioned in Lemma \ref{root}. We estimate \( {\mbox{\mathversion{bold}$u$}}_{2}\) in 
Sobolev--Slobodetski\u\i \ spaces. To estimate \( {\mbox{\mathversion{bold}$u$}}_{2}\) in \( H^{l, l/2}_{h}(Q_{\infty}) \),
we use the following norm.
\begin{eqnarray*}
\sum _{j\leq l}\int^{\infty}_{-\infty} \bigg\| \frac{\partial ^{j} \tilde{{\mbox{\mathversion{bold}$u$}}}_{2}}{\partial x^{j}}(\cdot ,\tau ) \bigg\|^{2}|\tau |^{l-j}{\rm d}\eta .
\end{eqnarray*}
Since 
\begin{eqnarray*}
\frac{\partial ^{j} \tilde{{\mbox{\mathversion{bold}$u$}}}_{2}}{\partial x^{j}} = \mu (\tau ) ^{j-1}\tilde{{\mbox{\mathversion{bold}$\Phi $}}}(\tau ){\rm e}^{\mu x},
\end{eqnarray*}
we have
\begin{eqnarray*}
\begin{aligned}
\bigg\| \frac{\partial ^{j} \tilde{{\mbox{\mathversion{bold}$u$}}}_{2}}{\partial x^{j}}(\cdot ,\tau ) \bigg\|^{2}
&=\int^{\infty}_{0} |\mu |^{2(j-1)}|\tilde{{\mbox{\mathversion{bold}$\Phi $}}}|^{2}|{\rm e}^{\mu x}|^{2}{\rm d}x \\
&= |\mu |^{2(j-1)}|\tilde{{\mbox{\mathversion{bold}$\Phi $}}}|^{2}\left( -\frac{1}{2\Re \mu }\right) .
\end{aligned}
\end{eqnarray*}
So we have
\begin{eqnarray*}
\int^{\infty}_{-\infty} \bigg\| \frac{\partial ^{j} \tilde{{\mbox{\mathversion{bold}$u$}}}_{2}}{\partial x^{j}}(\cdot ,\tau ) \bigg\|^{2}|\tau |^{l-j}{\rm d}\eta
=\int^{\infty}_{-\infty} |\tilde{{\mbox{\mathversion{bold}$\Phi $}}}(\tau )|^{2}|\mu (\tau ) |^{2(j-1)} \left( \frac{1}{2|\Re \mu (\tau )| }\right) |\tau |^{l-j}{\rm d}\eta .
\end{eqnarray*}
We divide the above integral into two parts, namely the part with \( |\eta | \geq \eta _{0}\) and \( |\eta |\leq \eta _{0}\), where 
\( \eta _{0} \) is defined in Lemma \ref{root}. From Lemma \ref{root}, in the domain \( |\eta |\geq \eta _{0}\), we have the asymptotic expansion
\begin{eqnarray*}
\left| \mu - \sqrt{\frac{\varepsilon}{2}}\left( -1+{\rm i}\right)|\eta |^{1/2} \right| \leq C,
\end{eqnarray*}
which implies, by taking \( \eta _{0}\) larger if necessary, \( |\frac{\tau }{|\mu |^{2}}|\leq C\). So we have
\begin{eqnarray*}
\int_{|\eta |\geq \eta _{0}}|\tilde{{\mbox{\mathversion{bold}$\Phi $}}}(\tau )|^{2}|\mu (\tau ) |^{2(j-1)} \left( \frac{1}{2|\Re \mu (\tau )| }\right) |\tau |^{l-j}{\rm d}\eta 
\leq C\int_{|\eta |\geq \eta _{0}} |\tilde{{\mbox{\mathversion{bold}$\Phi $}}}(\tau )|^{2}|\tau |^{l-3/2}{\rm d}\eta ,
\end{eqnarray*}
\begin{eqnarray*}
\int_{|\eta |\leq \eta _{0}}|\tilde{{\mbox{\mathversion{bold}$\Phi $}}}(\tau )|^{2}|\mu (\tau ) |^{2(j-1)} \left( \frac{1}{2|\Re \mu (\tau )| }\right) |\tau |^{l-j}{\rm d}\eta
\leq C\int^{\infty}_{-\infty}|\tilde{{\mbox{\mathversion{bold}$\Phi $}}}(\tau )|^{2}{\rm d}\eta.
\end{eqnarray*}
Combining these estimates we have
\begin{eqnarray*}
\begin{aligned}
||| {\mbox{\mathversion{bold}$u$}}_{2}||| ^{2}_{H^{l,l/2}_{h}(Q_{\infty})} &\leq C||| {\mbox{\mathversion{bold}$u $}}_{1x}|||^{2}_{H^{l-1,l/2-1/2}_{h}(Q_{\infty})} \\
&\leq C||| {\mbox{\mathversion{bold}$u $}}_{1}|||^{2}_{H^{l,l/2}_{h}(Q_{\infty})}.
\end{aligned}
\end{eqnarray*}
By choosing \( l=2k\) for an integer \( k\) and from Sobolev's imbedding theorem, we see that 
\begin{eqnarray}
\begin{aligned}
{\mbox{\mathversion{bold}$u$}}_{2}&\in H^{2k,k}_{h}(Q_{T})\hookrightarrow C\big( [0,T];H^{2k-2}(\mathbf{R}_{+})\big), \\
\frac{\partial ^{m}{\mbox{\mathversion{bold}$u$}}_{2}}{\partial t^{m}}&\in H^{2(k-m),k-m}_{h}(Q_{T}) \hookrightarrow C\big( [0,T];H^{2(k-m)-2}(\mathbf{R}_{+})\big),
\end{aligned}
\label{inc}
\end{eqnarray}
for \( m\leq k\).
We mentioned in section 3 that the given data can be taken as being smooth as desired while satisfying the necessary compatibility conditions, 
so the above construction of the solutions, estimates, and inclusions of function spaces  
imply that for an arbitrary integer \( l\), we have constructed a solution of (\ref{lp1}) such that 
\begin{eqnarray*}
{\mbox{\mathversion{bold}$u$}}\in \bigcap ^{l}_{j=0} C^{j}\big( [0,T]; H^{2(l-j)}(\mathbf{R}_{+})\big).
\end{eqnarray*}
To prove the uniqueness, we derive an energy estimate for the solution of (\ref{lp1}) directly. We calculate the following. 
\begin{eqnarray*}
\begin{aligned}
\frac{1}{2}\frac{{\rm d}}{{\rm d}t}\big\{ \| {\mbox{\mathversion{bold}$u$}}\| ^{2} + \alpha ^{2}\varepsilon ^{2}\| {\mbox{\mathversion{bold}$u$}}_{x}\| ^{2}\big\}
&= ({\mbox{\mathversion{bold}$u$}},{\mbox{\mathversion{bold}$u$}}_{t}) + \alpha ^{2}\varepsilon ^{2}({\mbox{\mathversion{bold}$u$}}_{x},{\mbox{\mathversion{bold}$u$}}_{tx})\\
&\leq -\alpha {\mbox{\mathversion{bold}$u$}}(0,t)\cdot {\mbox{\mathversion{bold}$u$}}_{xx}(0,t)
+\frac{\alpha \varepsilon }{2}|{\mbox{\mathversion{bold}$u$}}(0,t)|^{2}-\alpha ^{2}\varepsilon \| {\mbox{\mathversion{bold}$u$}}_{xx}\|^{2} \\
&\hspace*{1cm}+\frac{1}{2}\big( \| {\mbox{\mathversion{bold}$u$}}\| ^{2} + \alpha ^{2}\varepsilon ^{2}\| {\mbox{\mathversion{bold}$u$}}_{x}\| ^{2}\big) 
+ \| {\mbox{\mathversion{bold}$g$}} \|^{2},
\end{aligned}
\end{eqnarray*}
\begin{eqnarray*}
\begin{aligned}
\frac{1}{2}\frac{{\rm d}}{{\rm d}t}\big\{ \| {\mbox{\mathversion{bold}$u$}}_{x}\| ^{2} + \alpha ^{2}\varepsilon ^{2}\| {\mbox{\mathversion{bold}$u$}}_{xx}\| ^{2}\big\}
&= ({\mbox{\mathversion{bold}$u$}}_{x},{\mbox{\mathversion{bold}$u$}}_{xt}) + \alpha ^{2}\varepsilon ^{2}({\mbox{\mathversion{bold}$u$}}_{xx},{\mbox{\mathversion{bold}$u$}}_{xxt}) \\
&= \frac{\alpha }{2}|{\mbox{\mathversion{bold}$u$}}_{xx}(0,t)|^{2}-\alpha ^{2}\varepsilon \| {\mbox{\mathversion{bold}$u$}}_{xxx}\|^{2}
-( {\mbox{\mathversion{bold}$u$}}_{xx},{\mbox{\mathversion{bold}$g$}})-\alpha \varepsilon ( {\mbox{\mathversion{bold}$u$}}_{xxx},{\mbox{\mathversion{bold}$g$}}) \\
&\hspace*{1cm}-\alpha ^{2}\varepsilon {\mbox{\mathversion{bold}$u$}}_{xx}(0,t)\cdot {\mbox{\mathversion{bold}$u$}}_{xxx}(0,t) 
-\alpha \varepsilon {\mbox{\mathversion{bold}$u$}}_{xx}(0,t)\cdot {\mbox{\mathversion{bold}$g$}}(0,t) \\
&= \frac{\alpha }{2}|{\mbox{\mathversion{bold}$u$}}_{xx}(0,t)|^{2}-\alpha ^{2}\varepsilon \| {\mbox{\mathversion{bold}$u$}}_{xxx}\|^{2}
-( {\mbox{\mathversion{bold}$u$}}_{xx},{\mbox{\mathversion{bold}$g$}})-\alpha \varepsilon ( {\mbox{\mathversion{bold}$u$}}_{xxx},{\mbox{\mathversion{bold}$g$}}) \\
&\hspace*{1cm} -\alpha \varepsilon {\mbox{\mathversion{bold}$u$}}_{xx}(0,t)\cdot {\mbox{\mathversion{bold}$u$}}_{t}(0,t).
\end{aligned}
\end{eqnarray*}
On the other hand, we have from the equation,
\begin{eqnarray*}
\| {\mbox{\mathversion{bold}$u$}}_{t}+ \alpha \varepsilon {\mbox{\mathversion{bold}$u$}}_{tx}\|^{2}
= \| \alpha {\mbox{\mathversion{bold}$u$}}_{xxx}+{\mbox{\mathversion{bold}$g$}}\| ^{2}.
\end{eqnarray*}
Expanding the left-hand side gives
\begin{eqnarray*}
\begin{aligned}
\| {\mbox{\mathversion{bold}$u$}}_{t}+ \alpha \varepsilon {\mbox{\mathversion{bold}$u$}}_{tx}\|^{2}
&= \| {\mbox{\mathversion{bold}$u$}}_{t}\|^{2} + 2\alpha \varepsilon ({\mbox{\mathversion{bold}$u$}}_{t},{\mbox{\mathversion{bold}$u$}}_{tx})
+ \| {\mbox{\mathversion{bold}$u$}}_{tx}\| ^{2}\\
& = \| {\mbox{\mathversion{bold}$u$}}_{t}\|^{2} -\alpha \varepsilon |{\mbox{\mathversion{bold}$u$}}_{t}(0,t)|^{2} + \| {\mbox{\mathversion{bold}$u$}}_{tx}\| ^{2}.
\end{aligned}
\end{eqnarray*}
So we have 
\begin{eqnarray*}
-\alpha \varepsilon |{\mbox{\mathversion{bold}$u$}}_{t}(0,t)|^{2} \leq  \| \alpha {\mbox{\mathversion{bold}$u$}}_{xxx}+{\mbox{\mathversion{bold}$g$}}\| ^{2}.
\end{eqnarray*}
Utilizing the above estimate, we have for any positive \( \gamma \)
\begin{eqnarray*}
-\alpha \varepsilon |{\mbox{\mathversion{bold}$u$}}_{xx}(0,t)\cdot {\mbox{\mathversion{bold}$u$}}_{t}(0,t)| 
\leq -\alpha \gamma |{\mbox{\mathversion{bold}$u$}}_{xx}(0,t)|^{2} + \frac{5\alpha ^{2}\varepsilon }{18 \gamma }\| {\mbox{\mathversion{bold}$u$}}_{xxx}\|^{2}
+ C\| {\mbox{\mathversion{bold}$g$}}\| ^{2}.
\end{eqnarray*}
By choosing \( \frac{5}{18} < \gamma <\frac{1}{2} \), both \(  |{\mbox{\mathversion{bold}$u$}}_{xx}(0,t)|^{2}\) and \( \| {\mbox{\mathversion{bold}$u$}}_{xxx}\|^{2} \) 
can be dealt with in the estimates. Combining all the above estimates, we arrive at
\begin{eqnarray}
\| {\mbox{\mathversion{bold}$u$}}(t)\| ^{2}_{2} + \int ^{t}_{0} \| {\mbox{\mathversion{bold}$u$}}_{xx}(r)\| ^{2}_{1} {\rm d}r
\leq C\| {\mbox{\mathversion{bold}$u$}}_{0}\|^{2}_{2} + C\int ^{t}_{0} \| {\mbox{\mathversion{bold}$g$}}(r) \|^{2} {\rm d}r.
\label{fun}
\end{eqnarray}
By taking the \( t\) derivative, applying the same estimate as above, and estimating \( \| \partial ^{k}_{t}{\mbox{\mathversion{bold}$u$}}(\cdot, 0)\| \) as we did in the 
estimate of \( {\mbox{\mathversion{bold}$u$}}_{1}\), we have
\begin{eqnarray*}
\| \partial ^{k}_{t}{\mbox{\mathversion{bold}$u$}}(t) \| ^{2}_{2}+ \int^{t}_{0}\| \partial ^{k}_{t}{\mbox{\mathversion{bold}$u$}}_{xx}(r)\|^{2}_{1} {\rm d}r
\leq C\bigg( \| {\mbox{\mathversion{bold}$u$}}_{0}\|^{2}_{2+2k} + \sum ^{k-1}_{j=0}\| \partial ^{j}_{t}{\mbox{\mathversion{bold}$g$}}(t)\|^{2}_{2+2(k-1-j)}
+ \int ^{t}_{0}\| \partial ^{k}_{t}{\mbox{\mathversion{bold}$g$}}(r)\| ^{2} {\rm d}r\bigg).
\end{eqnarray*}
By using the equation to convert the time regularity into spacial regularity, we have for any \( k\) satisfying \( 0 \leq k\leq l\)
\begin{eqnarray*}
\begin{aligned}
&\sup _{0\leq t\leq T}\| \partial ^{k}_{t}{\mbox{\mathversion{bold}$u$}}(t)\| ^{2}_{2+2(l-k)} 
+ \int ^{T}_{0} \| \partial^{k}_{t} \mbox{\mathversion{bold}$u$}_{xx}(t)\|^{2}_{1+2(l-k)}{\rm d}t \\
&\hspace*{3cm}\leq C \bigg( \| {\mbox{\mathversion{bold}$u$}}_{0}\|^{2}_{2+2l} 
+ \sum ^{l-1}_{j=0}\sup _{0\leq t\leq T} \| \partial ^{j}_{t}{\mbox{\mathversion{bold}$g$}}(t)\|^{2}_{2+2(l-1-j)}
+ \int ^{T}_{0}\| \partial ^{l}_{t}{\mbox{\mathversion{bold}$g$}}(t)\| ^{2} {\rm d}t
\bigg).
\end{aligned}
\end{eqnarray*}
Up to this point, we have assumed that the given data are smooth. Through an approximation argument, we can relax the assumption on the data and prove the following.
\begin{lm}
For an arbitrary natural number \( l\), if \( {\mbox{\mathversion{bold}$u$}}_{0}\in H^{2+2l}(\mathbf{R}_{+})\),

\noindent $ \displaystyle {\mbox{\mathversion{bold}$g$}}\in \bigcap ^{l-1}_{j=0}C^{j}\big([0,T]; H^{2+2(l-1-j)}(\mathbf{R}_{+})\big) $, and
\( \partial ^{l}_{t}{\mbox{\mathversion{bold}$g$}} \in L^{2}(Q_{T})\)
 satisfy the compatibility conditions up to order \( l\),
there exists a unique solution \( {\mbox{\mathversion{bold}$u$}}\) to {\rm (\ref{lp1})} satisfying
\begin{eqnarray*}
\begin{aligned}
&\sum ^{l}_{k=0}\bigg( \sup _{0\leq t\leq T}\| \partial ^{k}_{t}{\mbox{\mathversion{bold}$u$}}(t)\| ^{2}_{2+2(l-k)} 
+ \int ^{T}_{0} \| \partial^{k}_{t} \mbox{\mathversion{bold}$u$}_{xx}(t)\|^{2}_{1+2(l-k)}{\rm d}t \bigg) \\
&\hspace*{3cm}\leq C \bigg( \| {\mbox{\mathversion{bold}$u$}}_{0}\|^{2}_{2+2l} + \sum ^{l-1}_{j=0}
 \sup _{0\leq t\leq T}\| \partial ^{j}_{t}{\mbox{\mathversion{bold}$g$}}(t)\|^{2}_{2+2(l-1-j)}
+ \int ^{T}_{0}\| \partial ^{l}_{t}{\mbox{\mathversion{bold}$g$}}(t)\| ^{2} {\rm d}t\bigg).
\end{aligned}
\end{eqnarray*}
\label{ex1}
\end{lm}
%


\section{Solving the Parabolic-Dispersive System}
\setcounter{equation}{0}

Now we construct the solution \( {\mbox{\mathversion{bold}$u$}}\) of (\ref{lpp1}) such that for a natural number \( l\)
\begin{eqnarray}
{\mbox{\mathversion{bold}$u$}}\in \bigcap ^{l}_{j=0}\left\{ C^{j}\big( [0,T]; H^{2+2(l-j)}(\mathbf{R}_{+})\big) 
\cap H^{j}\big(0,T;H^{3+2(l-j)}(\mathbf{R}_{+})\big) \right\},
\label{fs}
\end{eqnarray}
by iteration. For \( n\geq 1\), we define \( {\mbox{\mathversion{bold}$u$}}^{(n)}\) as the solution of the following 
system.
\begin{eqnarray*}
\left\{
\begin{array}{ll}
{\mbox{\mathversion{bold}$u$}}^{(n)}_{t}= \alpha {\mbox{\mathversion{bold}$u$}}^{(n)}_{xxx} + \alpha \varepsilon {\mbox{\mathversion{bold}$u$}}^{(n)}_{tx}
+{\rm A}({\mbox{\mathversion{bold}$w$}},\partial_{x}){\mbox{\mathversion{bold}$u$}}^{(n-1)}+ {\mbox{\mathversion{bold}$g$}}, & x>0, t>0, \\
{\mbox{\mathversion{bold}$u$}}^{(n)}(x,0)={\mbox{\mathversion{bold}$u$}}_{0}(x), & x>0, \\
{\mbox{\mathversion{bold}$u$}}^{(n)}_{x}(0,t)={\mbox{\mathversion{bold}$0$}}, & t>0.
\end{array}\right.
\end{eqnarray*}
\( {\mbox{\mathversion{bold}$u$}}^{(0)}\) must be defined in a specific way so that the compatibility conditions for each successive iteration will become satisfied. 
First, again by the approximation argument, we will assume that \( {\mbox{\mathversion{bold}$u$}}_{0}\) and \( {\mbox{\mathversion{bold}$g$}}\) satisfy the 
compatibility conditions for (\ref{lpp1}) up to order \( N\) and are smooth. We make the following notation.
\begin{eqnarray*}
Q({\mbox{\mathversion{bold}$v$}}):= \alpha {\mbox{\mathversion{bold}$v$}}_{xxx}
+ {\rm A}({\mbox{\mathversion{bold}$w$}},\partial_{x}){\mbox{\mathversion{bold}$v$}} +{\mbox{\mathversion{bold}$g$}}.
\end{eqnarray*}
We define \( {\mbox{\mathversion{bold}$u$}}^{(0)}\) as
\begin{eqnarray*}
{\mbox{\mathversion{bold}$u$}}^{(0)}(x,t):= {\mbox{\mathversion{bold}$u$}}_{0}(x)
+ \sum ^{N}_{j=1} \frac{t^j}{j!}\left( \frac{\partial ^{j}}{\partial t^{j}}Q({\mbox{\mathversion{bold}$v$}})\right) (x,0),
\end{eqnarray*}
where \( {\mbox{\mathversion{bold}$v$}}(x,0):= {\mbox{\mathversion{bold}$u$}}_{0}(x)\) and for \( k\geq 1\), 
\( {\mbox{\mathversion{bold}$\psi $}}_{k}(x):=\frac{\partial ^{k}{\mbox{\mathversion{bold}$v$}}}{\partial t^{k}}(x,0)\) is defined as the solution of the following 
linear ordinary differential equation, under the constraint that \( {\mbox{\mathversion{bold}$\psi $}}_{k}\) is integrable over \( \mathbf{R}_{+}\).
\begin{eqnarray*}
{\mbox{\mathversion{bold}$\psi $}}_{k}^{\prime }= -\frac{1}{\alpha \varepsilon }{\mbox{\mathversion{bold}$\psi $}}_{k}
+ \frac{1}{\alpha \varepsilon }\bigg( \alpha {\mbox{\mathversion{bold}$\psi $}}_{k-1}^{\prime \prime \prime}
 + \sum^{k-1}_{j=0}
\left(
\begin{array}{c}
k-1\\
j
\end{array}\right)
\big( \partial ^{j}_{t}
{\rm A}({\mbox{\mathversion{bold}$w$}}(\cdot ,0),\partial_{x})\big){\mbox{\mathversion{bold}$\psi $}}_{k-1-j} 
+ \partial ^{k-1}_{t}{\mbox{\mathversion{bold}$g$}}(\cdot ,0)\bigg).
\end{eqnarray*}
\( N\) is chosen to accommodate the necessary order of compatibility conditions and regularity. 
By defining \( {\mbox{\mathversion{bold}$u$}}^{(0)}\) as such, the compatibility conditions for each successive \( {\mbox{\mathversion{bold}$u$}}^{(n)}\) is 
automatically satisfied. 
From Lemma \ref{ex1}, \( \{ {\mbox{\mathversion{bold}$u$}}^{(n)}\}_{n\geq 1}\) is well-defined. We prove that 
\( \{ {\mbox{\mathversion{bold}$u$}}^{(n)}\} \) converges in the desired function space. From the way that we constructed \( {\mbox{\mathversion{bold}$u$}}^{(0)}\) we have
\begin{eqnarray*}
\begin{aligned}
\sum ^{l}_{k=0}\sup _{0\leq t\leq T}\| \partial ^{k}_{t}{\mbox{\mathversion{bold}$u$}}^{(0)}(t)\| ^{2}_{2+2(l-k)} 
&\leq C_{0}\left( \| {\mbox{\mathversion{bold}$u$}}_{0}\|^{2}_{2+2l+3N} + \sum ^{l-1}_{j=0}
\sup _{0\leq t\leq T}\| \partial ^{j}_{t}{\mbox{\mathversion{bold}$g$}}(t)\|^{2}_{2+2(l-1-j)+3N} \right). 
\end{aligned}
\end{eqnarray*}
Setting \( {\mbox{\mathversion{bold}$z$}}^{(n)}:= {\mbox{\mathversion{bold}$u$}}^{(n)}-{\mbox{\mathversion{bold}$u$}}^{(n-1)}\) we have
\begin{eqnarray*}
\left\{ 
\begin{array}{l}
{\mbox{\mathversion{bold}$z$}}^{(n)}_{t} = \alpha {\mbox{\mathversion{bold}$z$}}^{(n)}_{xxx}-\alpha \varepsilon {\mbox{\mathversion{bold}$z$}}^{(n)}_{tx}
+ {\rm A}({\mbox{\mathversion{bold}$w$}},\partial_{x}){\mbox{\mathversion{bold}$z$}}^{(n-1)}, \\
{\mbox{\mathversion{bold}$z$}}^{(n)}(x,0)={\mbox{\mathversion{bold}$0$}}, \\
{\mbox{\mathversion{bold}$z$}}^{(n)}_{x}(0,t)={\mbox{\mathversion{bold}$0$}}.
\end{array}\right.
\end{eqnarray*}
In the same way that we derived (\ref{fun}), we have
\begin{eqnarray*}
\begin{aligned}
\sup_{0\leq t\leq T}\| {\mbox{\mathversion{bold}$z$}}^{(n)}(t)\|^{2}_{2} + \int ^{T}_{0}\| {\mbox{\mathversion{bold}$z$}}^{(n)}_{xx}(t)\| ^{2}_{1}{\rm d}t
&\leq C\int ^{T}_{0}\| {\mbox{\mathversion{bold}$z$}}^{(n-1)}_{xx}(t)\| ^{2}{\rm d}t \\
&\leq \frac{(CT)^{n-1}}{(n-1)!}M.
\end{aligned}
\end{eqnarray*}
The above estimate proves that \( {\mbox{\mathversion{bold}$u$}}^{(n)}\) converges in \( C\big([0,T];H^{2}(\mathbf{R}_{+})\big) \cap L^{2}\big(0,T;H^{3}(\mathbf{R}_{+})\big) \).
Since \( \partial ^{k}_{t}{\mbox{\mathversion{bold}$z$}}^{(n)}(x,0)={\mbox{\mathversion{bold}$0$}}\), we can prove in the same way as above that
\( \partial ^{k}_{t}{\mbox{\mathversion{bold}$u$}}^{(n)} \) converges in \( C\big([0,T];H^{2}(\mathbf{R}_{+})\big) \cap L^{2}\big(0,T;H^{3}(\mathbf{R}_{+})\big) \).
Using the equation we can prove that for \( 0\leq k \leq l \), \( \partial ^{k}_{t}{\mbox{\mathversion{bold}$u$}}^{(n)}\) converges in 
\( C\left([0,T];H^{2+2(l-k)}(\mathbf{R}_{+})\right) \cap L^{2}\big( 0,T; H^{3+2(l-k)}(\mathbf{R}_{+})\big) \). 
Thus, for an arbitrary \( l\), we have constructed a solution of (\ref{lpp1}) satisfying the condition (\ref{fs}).

Now we consider the limit \( \varepsilon \rightarrow 0\). For this, we derive an estimate of the solution that is uniform in \( \varepsilon \). 
The energy form we use is the same as the estimate we obtained before, but we use the elliptic term to make the estimate uniform in \( \varepsilon \).
We are still assuming that the given data are smooth as necessary.
We estimate as follows.
\begin{eqnarray*}
\begin{aligned}
\frac{1}{2}\frac{{\rm d}}{{\rm d}t}\big\{ \|{\mbox{\mathversion{bold}$u$}}\| ^{2} + \alpha ^{2}\varepsilon ^{2}\|{\mbox{\mathversion{bold}$u$}}_{x}\|^{2}\big\}
&\leq -\alpha{\mbox{\mathversion{bold}$u$}}(0,t)\cdot {\mbox{\mathversion{bold}$u$}}_{xx}(0,t) 
+\frac{\alpha \varepsilon}{2}\frac{{\rm d}}{{\rm d}t}|{\mbox{\mathversion{bold}$u$}}(0,t)|^{2}
-\frac{\delta}{2} \|{\mbox{\mathversion{bold}$u$}}_{x}\|^{2}-\alpha ^{2}\varepsilon \|{\mbox{\mathversion{bold}$u$}}_{xx}\|^{2}	\\
&\hspace*{0.5cm} + C\|{\mbox{\mathversion{bold}$u$}}\|^{2} +\alpha ^{2}\varepsilon ^{2}\|{\mbox{\mathversion{bold}$u$}}_{x}\|^{2}
+\varepsilon ^{2}\|{\rm A}({\mbox{\mathversion{bold}$w$}},\partial_{x}){\mbox{\mathversion{bold}$u$}}\|^{2}.
\end{aligned}
\end{eqnarray*}
We choose \( \varepsilon _{1}>0\) such that \( \varepsilon _{1}\|{\rm A}_{0}({\mbox{\mathversion{bold}$w$}})\|^{2}_{L^{\infty}(0,T;L^{\infty}(\mathbf{R}_{+}))}
\leq \frac{\alpha ^{2}}{2}\).
Then, for \( 0<\varepsilon \leq \varepsilon _{1}\) we have 
\begin{eqnarray*}
\begin{aligned}
\frac{1}{2}\frac{{\rm d}}{{\rm d}t}\big\{ \|{\mbox{\mathversion{bold}$u$}}\| ^{2} + \alpha ^{2}\varepsilon ^{2}\|{\mbox{\mathversion{bold}$u$}}_{x}\|^{2}\big\}
&\leq -\alpha {\mbox{\mathversion{bold}$u$}}(0,t)\cdot {\mbox{\mathversion{bold}$u$}}_{xx}(0,t)
+\frac{\alpha \varepsilon }{2}\frac{{\rm d}}{{\rm d}t}|{\mbox{\mathversion{bold}$u$}}(0,t)|^{2}-\frac{\delta}{2}\| {\mbox{\mathversion{bold}$u$}}_{x}\|^{2}
- \frac{\alpha ^{2}\varepsilon }{2}\| {\mbox{\mathversion{bold}$u$}}_{xx}\|^{2}\\[0.3cm]
&\hspace*{1cm} +C\|{\mbox{\mathversion{bold}$u$}}\|_{1}^{2}
 + \|{\mbox{\mathversion{bold}$g$}}\|^{2}.
\end{aligned}
\end{eqnarray*}
Next, we have
\begin{eqnarray*}
\begin{aligned}
\frac{1}{2}\frac{{\rm d}}{{\rm d}t}\left\{ \|{\mbox{\mathversion{bold}$u$}}_{x}\|^{2} + \alpha ^{2}\varepsilon ^{2}\| {\mbox{\mathversion{bold}$u$}}_{xx}\|^{2}\right\}
&\leq \frac{\alpha }{6}|{\mbox{\mathversion{bold}$u$}}_{xx}(0,t)|^{2} -\frac{\delta}{2} \|{\mbox{\mathversion{bold}$u$}}_{xx}\|^{2}
-\frac{\alpha ^{2}\varepsilon }{12}\|{\mbox{\mathversion{bold}$u$}}_{xxx}\|^{2}\\
&\hspace*{0.4cm}+\frac{3\alpha \varepsilon }{4}\|{\rm A}_{0}({\mbox{\mathversion{bold}$w$}})\|^{2}_{L^{\infty}(0,T;L^{\infty}(\mathbf{R}_{+}))}
\|{\mbox{\mathversion{bold}$u$}}_{x}\| \|{\mbox{\mathversion{bold}$u$}}_{xxx}\| \\
&\hspace*{0.4cm}-\alpha \varepsilon ({\mbox{\mathversion{bold}$u$}}_{xxx},{\rm A}({\mbox{\mathversion{bold}$w$}},\partial_{x}){\mbox{\mathversion{bold}$u$}}) 
+C\big( \|{\mbox{\mathversion{bold}$u$}}\|_{1}^{2}+ \|{\mbox{\mathversion{bold}$g$}}\|^{2}\big),
\end{aligned}
\end{eqnarray*}
where we have used the interpolation inequality \( \| \mbox{\mathversion{bold}$u$}_{xx}\|^{2}\leq C\| \mbox{\mathversion{bold}$u$}_{x}\| \|\mbox{\mathversion{bold}$u$}_{xxx}\|\) . 
Now we choose \( \varepsilon _{2}>0\) so that\\ 
\( \frac{3}{8}\varepsilon _{2}\|{\rm A}_{0}({\mbox{\mathversion{bold}$w$}})\|^{2}_{L^{\infty}(0,T;L^{\infty}(\mathbf{R}_{+}))}
\leq \frac{\alpha }{24}\)
and \( 24\varepsilon _{2}\| {\rm A}_{0}({\mbox{\mathversion{bold}$w$}})\|^{2}_{L^{\infty}(0,T;L^{\infty}(\mathbf{R}_{+}))}\leq \frac{\delta }{4}\).
Then, for \( 0\leq \varepsilon \leq \varepsilon _{2}\) we have
\begin{eqnarray*}
\frac{1}{2}\frac{{\rm d}}{{\rm d}t}\left\{ \|{\mbox{\mathversion{bold}$u$}}_{x}\|^{2} + \alpha ^{2}\varepsilon ^{2}\| {\mbox{\mathversion{bold}$u$}}_{xx}\|^{2}\right\}
\leq \frac{\alpha }{6}|{\mbox{\mathversion{bold}$u$}}_{xx}(0,t)|^{2}-\frac{\delta}{4}\|{\mbox{\mathversion{bold}$u$}}_{xx}\|^{2}
-\frac{\alpha ^{2}\varepsilon}{24}\|{\mbox{\mathversion{bold}$u$}}_{xxx}\|^{2}+ C\big( \|{\mbox{\mathversion{bold}$u$}}_{x}\|^{2}+ \|{\mbox{\mathversion{bold}$g$}}\|^{2}\big).
\end{eqnarray*}
Finally we estimate
\begin{eqnarray*}
\frac{1}{2}\frac{{\rm d}}{{\rm d}t}\|{\mbox{\mathversion{bold}$u$}}_{xx}\|^{2}
\leq \frac{\alpha }{2}|{\mbox{\mathversion{bold}$u$}}_{xxx}(0,t)|^{2} -\frac{\varepsilon }{2}\|{\mbox{\mathversion{bold}$u$}}_{tx}\|^{2}
-\frac{\delta }{4}\| {\mbox{\mathversion{bold}$u$}}_{xxx}\|^{2} + C\|{\mbox{\mathversion{bold}$u$}}\|_{2}^{2} + \|{\mbox{\mathversion{bold}$g$}}_{x}\|^{2}.
\end{eqnarray*}
In each estimate, the constant \( C\) is independent of \( \varepsilon \in (0, \varepsilon _{0}]\), where \( \varepsilon _{0}:=\min \{\varepsilon _{1},\varepsilon _{2}\} \).
Combining all the estimates yields for \( 0\leq \varepsilon \leq \varepsilon _{0}\), 
\begin{eqnarray*}
\sup_{0\leq t\leq T} \| {\mbox{\mathversion{bold}$u$}}(t)\|^{2}_{2} + \delta \int ^{T}_{0} \big( \| {\mbox{\mathversion{bold}$u$}}_{xxx}(t)\|^{2}
+ \varepsilon \|{\mbox{\mathversion{bold}$u$}}_{tx}(t)\|^{2} + \alpha |{\mbox{\mathversion{bold}$u$}}_{xx}(0,t)|^{2} + \alpha |{\mbox{\mathversion{bold}$u$}}_{xxx}(0,t)|^{2}
\big){\rm d}t \\
\leq C\bigg( \|{\mbox{\mathversion{bold}$u$}}_{0}\|^{2}_{2} + \int ^{T}_{0} \| {\mbox{\mathversion{bold}$g$}}(t)\|^{2}_{1}{\rm d}t \bigg).
\end{eqnarray*}
Now we take the derivative of the equation \( m\) times \( (1\leq m\leq l)\) with respect to \( t\) and set 
\( {\mbox{\mathversion{bold}$v$}}_{m}:=\partial ^{m}_{t}{\mbox{\mathversion{bold}$u$}} \). Then,
\( {\mbox{\mathversion{bold}$v$}}_{m}\) satisfies
\begin{eqnarray*}
\left\{
\begin{array}{l}
{\mbox{\mathversion{bold}$v$}}_{mt} = \alpha {\mbox{\mathversion{bold}$v$}}_{mxxx} - \alpha \varepsilon 
{\rm A}({\mbox{\mathversion{bold}$w$}},\partial_{x}){\mbox{\mathversion{bold}$v$}}_{m} + \partial ^{m}_{t}{\mbox{\mathversion{bold}$g$}} + {\mbox{\mathversion{bold}$F$}}_{m},\\
{\mbox{\mathversion{bold}$v$}}_{m}(x,0)={\mbox{\mathversion{bold}$\phi $}}_{m}(x), \\
{\mbox{\mathversion{bold}$v$}}_{mx}(0,t)= {\mbox{\mathversion{bold}$0$}},
\end{array}\right.
\end{eqnarray*}
where $ \displaystyle {\mbox{\mathversion{bold}$F$}}_{m}=\sum ^{m-1}_{j=0} 
\left( 
\begin{array}{c}
m-1\\
j
\end{array}\right)
\partial^{m-1-j}_{t}\left({\rm A}({\mbox{\mathversion{bold}$w$}},\partial_{x})\right)
{\mbox{\mathversion{bold}$v$}}_{j} $. 
We derive the uniform estimate by induction on \( m\). The case \( m=0\) was just derived, so suppose for \( 0\leq j\leq m-1\) we have
\begin{eqnarray*}
\begin{aligned}
&\sup_{0\leq t\leq T} \| {\mbox{\mathversion{bold}$v$}}_{j}(t)\|^{2}_{2} + \delta \int^{T}_{0}\| {\mbox{\mathversion{bold}$v$}}_{jxxx}(t)\|^{2}{\rm d}t \\
&\hspace*{0.5cm}\leq C\left\{ \| {\mbox{\mathversion{bold}$u$}}_{0}\|^{2}_{2+3j} 
+ \sup_{0\leq t\leq T} \bigg(\|\partial ^{j-1}_{t}{\mbox{\mathversion{bold}$g$}}(\cdot ,t)\|^{2}_{2}\right.
+ \sum^{j-2}_{k=0}\|\partial ^{k}_{t}{\mbox{\mathversion{bold}$g$}}(\cdot, t)\|^{2}_{2+3(m-2-k)} \bigg)
\left. + \int ^{T}_{0}\| \partial ^{j}_{t}{\mbox{\mathversion{bold}$g$}}(t)\|_{1} {\rm d}t \right\}
\end{aligned}
\end{eqnarray*}
with \( C\) independent of \( \varepsilon \). 
Estimating in the same way as before, we have
\begin{eqnarray*}
\| {\mbox{\mathversion{bold}$v$}}_{m}(t)\|^{2}_{2} + \int ^{T}_{0}\| {\mbox{\mathversion{bold}$v$}}_{mxxx}(t)\|^{2}{\rm d}t
\leq C\bigg( \| {\mbox{\mathversion{bold}$\phi $}}_{m}\|^{2}_{2}+ \int ^{T}_{0}\big( \| \partial ^{m}_{t}{\mbox{\mathversion{bold}$g$}}(t)\|^{2}_{1}
+ \| {\mbox{\mathversion{bold}$F$}}_{m}(t)\|^{2}_{1}\big){\rm d}t\bigg).
\end{eqnarray*}
Now we estimate the right-hand side.
\begin{eqnarray*}
\| {\mbox{\mathversion{bold}$F$}}_{m}(t)\|_{1}^{2}\leq C \sum^{m-1}_{j=0}\|{\mbox{\mathversion{bold}$v$}}_{jxx}(t)\|^{2}_{1},
\end{eqnarray*}
where \( C\) depends on the norm of \({\mbox{\mathversion{bold}$w$}}\) in $ \displaystyle \bigcap ^{m-1}_{j=0}W^{j,\infty}\big( 0,T;H^{1}(\mathbf{R}_{+})\big) $.
The expression for \( {\mbox{\mathversion{bold}$\phi $}}_{m}\) and its derivatives are
\begin{eqnarray*}
\partial ^{k}_{x}{\mbox{\mathversion{bold}$\phi $}}_{m}(x)=-\frac{1}{\alpha \varepsilon }\int ^{\infty}_{x}{\rm e}^{-\frac{1}{\alpha \varepsilon }(x-y)}
\partial ^{k}_{y}\big\{\alpha {\mbox{\mathversion{bold}$\phi $}}^{\prime \prime \prime}_{m-1}(y)+ {\mbox{\mathversion{bold}$F$}}_{m-1}(y,0) 
+ \partial ^{m-1}_{t}{\mbox{\mathversion{bold}$g$}}(y,0)\big\}{\rm d}y.
\end{eqnarray*}
Through direct calculation, we see that 
\begin{eqnarray*}
\left\| -\frac{1}{\alpha \varepsilon }\int ^{\infty}_{\cdot } {\rm e}^{-\frac{1}{\alpha \varepsilon }(\cdot -y)}{\mbox{\mathversion{bold}$\Phi $}}(y){\rm d}y\right\|
\leq \|{\mbox{\mathversion{bold}$\Phi $}}\|.
\end{eqnarray*}
Thus, we can prove that 
\begin{eqnarray*}
\| {\mbox{\mathversion{bold}$\phi $}}_{m}\|_{2} \leq C\bigg( \| {\mbox{\mathversion{bold}$u$}}_{0}\|_{2+3m} + \| \partial ^{m-1}_{t}{\mbox{\mathversion{bold}$g$}}(\cdot ,0)\|_{2}
+ \sum ^{m-2}_{j=0}\|\partial ^{j}_{t}{\mbox{\mathversion{bold}$g$}}(\cdot ,0)\|_{2+3(m-2-j)}\bigg).
\end{eqnarray*}
Here, \( C\) depends on the norm of \( {\mbox{\mathversion{bold}$w$}}\) in $ \displaystyle \bigcap ^{m-1}_{j=0}C^{j}\big([0,T]; H^{2+3(m-1-j)}(\mathbf{R}_{+})\big) $ and 
the norm of \( \partial ^{m}_{t}{\mbox{\mathversion{bold}$w$}}\) in \( L^{\infty}\big( 0,T; H^{1}(\mathbf{R}_{+})\big) \), but is independent of \( \varepsilon \).
Combining these estimates and using the equation yields
\begin{eqnarray}
\begin{aligned}
&\sum ^{l}_{j=0}\bigg\{  \sup_{0\leq t\leq T} \| \partial ^{j}_{t}{\mbox{\mathversion{bold}$u$}}(t)\|^{2}_{2+3(l-j)}
+ \int ^{T}_{0} \| \partial ^{j}_{t}{\mbox{\mathversion{bold}$u$}}(t)\|^{2}_{3+3(l-j)}{\rm d}t \bigg\}\\
&\hspace*{1cm}\leq C\bigg\{ \| {\mbox{\mathversion{bold}$u$}}_{0}\|^{2}_{2+3l} 
+ \sum ^{l-1}_{j=0} \bigg( \sup_{0\leq t\leq T} \|\partial ^{j}_{t}{\mbox{\mathversion{bold}$g$}}(t)\|^{2}_{2+3(l-1-j)} 
+ \int ^{T}_{0} \| \partial ^{j}_{t}{\mbox{\mathversion{bold}$g$}}(t)\|^{2}_{3+3(l-1-j)}{\rm d}t \bigg) \\
&\hspace*{8cm} + \int^{T}_{0}\|\partial ^{l}_{t}{\mbox{\mathversion{bold}$g$}}(t)\|_{1}^{2}{\rm d}t s\bigg\}. 	
\end{aligned}
\label{fin1}	
\end{eqnarray}
Again, we emphasize that \( C\) is independent of \( \varepsilon \). 
Now we denote the solution of (\ref{lpp1}) as \( {\mbox{\mathversion{bold}$u$}}^{\varepsilon }\) to emphasize that the solution depends on \( \varepsilon \).
We also recall that \( {\mbox{\mathversion{bold}$g$}}\) was a correction of \( {\mbox{\mathversion{bold}$f$}}\) which depends on \( \varepsilon \), so 
we denote it as \( {\mbox{\mathversion{bold}$g$}}^{\varepsilon }\). 
To make the following arguments more simple, we will make the following notations. Recall that
\( X^{l}\) is the function space that we are constructing the solution in, specifically,
\begin{eqnarray}
X^{l}:= \bigcap ^{l}_{j=0} \bigg( C^{j}\big( [0,T];H^{2+3(l-j)}(\mathbf{R}_{+})\big) \cap H^{j}\big( 0,T;H^{3+3(l-j)}(\mathbf{R}_{+})\big) \bigg).
\label{X}
\end{eqnarray}
\( Y^{l}\) is the function space that \( {\mbox{\mathversion{bold}$f$}}\) will be required to belong in, and is defined by 
\begin{eqnarray}
Y^{l}:= \bigg\{ f; \  f\in \bigcap ^{l-1}_{j=0} C^{j}\big( [0,T];H^{2+3(l-1-j)}(\mathbf{R}_{+})\big) , \ \partial ^{l}_{t}f \in L^{2}\big( 0,T; H^{1}(\mathbf{R}_{+})\big) \bigg\}.
\label{Y}
\end{eqnarray}
Finally, \( Z^{l}\) is the function space that \( {\mbox{\mathversion{bold}$w$}}\) will belong in and is defined as
\begin{eqnarray}
Z^{l}:= \bigg\{ w; \ w\in \bigcap ^{l-1}_{j=0}C^{j}\big( [0,T];H^{2+3(l-1-j)}(\mathbf{R}_{+})\big), \ \partial ^{l}_{t}w\in L^{\infty}\big (0,T;H^{1}(\mathbf{R}_{+})\big) \bigg\}.
\label{Z}
\end{eqnarray}
We assume that \( {\mbox{\mathversion{bold}$u$}}_{0}\in H^{2+3N}(\mathbf{R}_{+})\), \ \( {\mbox{\mathversion{bold}$g$}}^{\varepsilon }\in Y^{N}\), 
\( {\mbox{\mathversion{bold}$w$}}\in Z^{N} \) for \( N>l+1\), and \( {\mbox{\mathversion{bold}$g$}}^{\varepsilon }\rightarrow {\mbox{\mathversion{bold}$f$}}\) in 
\( Y^{l+1}\). Thus we know the existence of a unique solution \( {\mbox{\mathversion{bold}$u$}}^{\varepsilon }\in X^{l+1}\) of (\ref{lpp1}) with a uniform bound in \( X^{l+1}\).  
For \( 0 < \varepsilon <\varepsilon ^{\prime} \leq \varepsilon _{0}\), we set 
\( {\mbox{\mathversion{bold}$z$}}:={\mbox{\mathversion{bold}$u$}}^{\varepsilon^{\prime} }-{\mbox{\mathversion{bold}$u$}}^{\varepsilon}\). 
\( {\mbox{\mathversion{bold}$z$}}\) satisfies
\begin{eqnarray*}
\left\{ 
\begin{array}{l}
{\mbox{\mathversion{bold}$z$}}_{t}=\alpha {\mbox{\mathversion{bold}$z$}}_{xxx} - \alpha \varepsilon ^{\prime}{\mbox{\mathversion{bold}$z$}}_{tx}
+ {\rm A}({\mbox{\mathversion{bold}$w$}},\partial_{x}){\mbox{\mathversion{bold}$z$}} 
-\alpha (\varepsilon ^{\prime} -\varepsilon ){\mbox{\mathversion{bold}$u$}}^{\varepsilon }_{tx} 
+ \mbox{\mathversion{bold}$g$}^{\varepsilon ^{\prime}}-\mbox{\mathversion{bold}$g$}^{\varepsilon }, \\
{\mbox{\mathversion{bold}$z$}}(x,0)={\mbox{\mathversion{bold}$0$}},\\
{\mbox{\mathversion{bold}$z$}}_{x}(0,t)={\mbox{\mathversion{bold}$0$}}.
\end{array}\right.
\end{eqnarray*}
From (\ref{fin1}), we have
\begin{eqnarray*}
\begin{aligned}
\| {\mbox{\mathversion{bold}$z$}}\|^{2}_{X^{l}}
&\leq C(\varepsilon^{\prime}+\varepsilon )^{2}\bigg\{  
 \sum ^{l-1}_{j=0} \bigg( \sup_{0\leq t\leq T} \|\partial ^{j+1}_{t}{\mbox{\mathversion{bold}$u$}}^{\varepsilon }_{x}(t)\|^{2}_{2+3(l-1-j)} 
+ \int ^{T}_{0} \| \partial ^{j+1}_{t}{\mbox{\mathversion{bold}$u$}}^{\varepsilon }_{x}(t)\|^{2}_{3+3(l-1-j)}{\rm d}t \bigg) \\
&\hspace*{3cm}+ \int^{T}_{0}\|\partial ^{l+1}_{t}{\mbox{\mathversion{bold}$u$}}^{\varepsilon }_{x}(t)\|_{1}^{2}{\rm d}t  
\bigg\} + \int ^{T}_{0} \| (\mbox{\mathversion{bold}$g$}^{\varepsilon ^{\prime}}-\mbox{\mathversion{bold}$g$}^{\varepsilon })(t)\|^{2}_{Y^{l+1}}{\rm d}t \\
&\leq C(\varepsilon ^{\prime}
+\varepsilon )^{2}+ \int ^{T}_{0} \| (\mbox{\mathversion{bold}$g$}^{\varepsilon ^{\prime}}-\mbox{\mathversion{bold}$g$}^{\varepsilon })(t)\|^{2}_{Y^{l+1}}{\rm d}t.
\end{aligned}
\end{eqnarray*}
Thus we see that there exists a \( {\mbox{\mathversion{bold}$u$}}\) such that \( {\mbox{\mathversion{bold}$u$}}^{\varepsilon }\rightarrow {\mbox{\mathversion{bold}$u$}}\) 
in \( X^{l}\), and \( {\mbox{\mathversion{bold}$u$}}\) is a solution of (\ref{lpd1}).	
We derive an energy estimate for \( {\mbox{\mathversion{bold}$u$}}\) to prove the uniqueness of the solution. 
Through a standard energy estimate, we obtain the following estimates.
\begin{eqnarray*}
\begin{aligned}
\frac{1}{2}\frac{{\rm d}}{{\rm d}t}\|{\mbox{\mathversion{bold}$u$}}\|^{2} &\leq -\alpha {\mbox{\mathversion{bold}$u$}}(0,t)\cdot {\mbox{\mathversion{bold}$u$}}_{xx}(0,t)
-\frac{\delta }{2}\|{\mbox{\mathversion{bold}$u$}}_{xx}\|^{2} + C\| {\mbox{\mathversion{bold}$u$}}\|^{2} + \|{\mbox{\mathversion{bold}$f$}}\|^{2},\\
\frac{1}{2}\frac{{\rm d}}{{\rm d}t}\|{\mbox{\mathversion{bold}$u$}}_{x}\|^{2} &\leq \frac{\alpha }{2}|{\mbox{\mathversion{bold}$u$}}_{xx}(0,t)|^{2}
-\frac{\delta }{2}\|{\mbox{\mathversion{bold}$u$}}_{xx}\|^{2} + C\big(\|{\mbox{\mathversion{bold}$u$}}\|^{2}_{1} + \| {\mbox{\mathversion{bold}$f$}}\|^{2}\big),\\
\frac{1}{2}\frac{{\rm d}}{{\rm d}t}\|{\mbox{\mathversion{bold}$u$}}_{xx}\|^{2} &\leq \frac{\alpha }{2}|{\mbox{\mathversion{bold}$u$}}_{xxx}(0,t)|^{2}
-\frac{\delta }{2}\|{\mbox{\mathversion{bold}$u$}}_{xxx}\|^{2} + C\big(\|{\mbox{\mathversion{bold}$u$}}\|^{2}_{2} + \| {\mbox{\mathversion{bold}$f$}}_{x}\|^{2}\big).
\end{aligned}
\end{eqnarray*}
So combining these estimates, we have
\begin{eqnarray*}
\sup_{0\leq t\leq T}\| {\mbox{\mathversion{bold}$u$}}(t)\|^{2}_{2} + \int ^{T}_{0}\|{\mbox{\mathversion{bold}$u$}}_{x}(t)\|^{2}_{2}{\rm d}t
\leq C\bigg( \| {\mbox{\mathversion{bold}$u$}}_{0}\|^{2}_{2} + \int ^{T}_{0} \| {\mbox{\mathversion{bold}$f$}}(t)\|^{2}_{1}{\rm d}t \bigg).
\end{eqnarray*}
As before, taking the derivative with respect to \( t\) in the equation, applying the above estimate, and converting the regularity in \( t\) into \( x\) via the equation, we have
\begin{eqnarray}
\| {\mbox{\mathversion{bold}$u$}}\|_{X^{l}}\leq C\big( \|{\mbox{\mathversion{bold}$u$}}_{0}\|_{2+3l} + \| {\mbox{\mathversion{bold}$f$}}\|_{Y^{l}}\big).
\label{finn1}
\end{eqnarray}
Here, \( C\) depends on \( \|{\mbox{\mathversion{bold}$w$}}\|_{Z^{l}}\), \( T\), and \( \delta \).

As in Lemma \ref{ex1}, we can relax the condition on the given data by taking approximating series 
\( \mbox{\mathversion{bold}$u$}_{0n}\subset H^{2+3N}\),
\( \{ {\mbox{\mathversion{bold}$f$}}_{n}\}_{n\geq 1} \subset Y^{N}\), and \( \{ {\mbox{\mathversion{bold}$w$}}_{n}\}_{n\geq 1} \subset Z^{N}\) which 
\( \mbox{\mathversion{bold}$u$}_{0n}\rightarrow \mbox{\mathversion{bold}$u$}_{0}\) in \( H^{2+3l}(\mathbf{R}_{+})\),
\( {\mbox{\mathversion{bold}$f$}}_{n}\rightarrow {\mbox{\mathversion{bold}$f$}}\) in \( Y^{l}\), and \( {\mbox{\mathversion{bold}$w$}}_{n} \rightarrow {\mbox{\mathversion{bold}$w$}}\)
in \( Z^{l}\). Applying (\ref{finn1}), and passing to the limit, we arrive at our main theorem. For notations of function spaces, see 
(\ref{X}), (\ref{Y}), and (\ref{Z}).

\begin{Th}
For any \( T>0\) and an arbitrary non-negative integer \( l\), if \( {\mbox{\mathversion{bold}$u$}}_{0}\in H^{2+3l}(\mathbf{R}_{+})\), 
\( {\mbox{\mathversion{bold}$f$}}\in Y^{l}\), and \( {\mbox{\mathversion{bold}$w$}}\in Z^{l}\) satisfy the compatibility conditions up to order \( l\), 
a unique solution \( {\mbox{\mathversion{bold}$u$}}\) of {\rm (\ref{lpd1})} exists such that \( {\mbox{\mathversion{bold}$u$}}\in X^{l}\).
\end{Th}
%
%
%
%
%
%
%
%
%
\section{The Case \( \alpha >0\)}
\setcounter{equation}{0}
The case \( \alpha >0\) can be treated by a standard argument. We start by considering the following regularized problem.
\begin{eqnarray}
\left\{
\begin{array}{ll}
\mbox{\mathversion{bold}$u$}_{t}=-\varepsilon\mbox{\mathversion{bold}$u$}_{xxxx}+\mbox{\mathversion{bold}$g$}, & x>0,t>0,\\
\mbox{\mathversion{bold}$u$}(x,0)=\mbox{\mathversion{bold}$u$}_{0}(x), & x>0,\\
\mbox{\mathversion{bold}$u$}(0,t)=\mbox{\mathversion{bold}$e$},& t>0,\\
\mbox{\mathversion{bold}$u$}_{x}(0,t)=\mbox{\mathversion{bold}$0$}, & t>0.
\end{array}\right.
\label{nr}
\end{eqnarray}
As before, we explicitly construct the solution of (\ref{nr}) in the form \( \mbox{\mathversion{bold}$u$}=\mbox{\mathversion{bold}$u$}^{1}+\mbox{\mathversion{bold}$u$}^{2}\). 
Where \( \mbox{\mathversion{bold}$u$}^{1}\) is defined as the solution of 
\begin{eqnarray*}
\left\{
\begin{array}{ll}
\mbox{\mathversion{bold}$u$}^{1}_{t}=-\varepsilon \mbox{\mathversion{bold}$u$}^{1}_{xxxx} + \mbox{\mathversion{bold}$G$}, & x\in \mathbf{R}, t>0,\\
\mbox{\mathversion{bold}$u$}^{1}(x,0)=\mbox{\mathversion{bold}$U$}_{0}(x),& x\in \mathbf{R},
\end{array}\right.
\end{eqnarray*}
and \( \mbox{\mathversion{bold}$u$}^{2}\) is defined as the solution of 
\begin{eqnarray*}
\left\{
\begin{array}{ll}
\mbox{\mathversion{bold}$u$}^{2}_{t}=-\varepsilon \mbox{\mathversion{bold}$u$}^{2}_{xxxx},& x>0,t>0,\\
\mbox{\mathversion{bold}$u$}^{2}(x,0)=\mbox{\mathversion{bold}$0$}, & x>0,\\
\mbox{\mathversion{bold}$u$}^{2}(0,t)=\mbox{\mathversion{bold}$e$}-\mbox{\mathversion{bold}$u$}^{1}(0,t), & t>0,\\
\mbox{\mathversion{bold}$u$}^{2}_{x}(0,t)= -\mbox{\mathversion{bold}$u$}^{1}_{x}(0,t), & t>0.
\end{array}\right.
\end{eqnarray*}
The solutions can be constructed using Fourier transform and Laplace transform as in the case \( \alpha <0\). 
We note that in estimating \( \mbox{\mathversion{bold}$u$}^{2}\), we slightly modifiy the Sobolev--Slobodetski\u\i \ space to 
fit our fourth order parabolic system. For an integer \( m\), we define the space \( H^{m,m/4}_{h}(Q_{T})\) analogous to 
\( H^{m,m/2}_{h}(Q_{T})\), and we use the case \( m=4l\) and the norm
\begin{eqnarray*}
\| \mbox{\mathversion{bold}$u$}\|_{  H^{4l,l}_{h}(Q_{T})}^{2}=
\sum _{j\leq l}\int ^{\infty}_{-\infty}\left\| \frac{\partial ^{4j}\tilde{\mbox{\mathversion{bold}$u$}}}{\partial x^{4j}}(\cdot ,\tau )\right\| ^{2}|\tau | ^{l-j}{\rm d}\eta.
\end{eqnarray*}
Then we construct the solution 
$\displaystyle\mbox{\mathversion{bold}$u$} \in \bigcap ^{l}_{j=0}C^{j}\big( [0,T];H^{2+4(l-j)}(\mathbf{R}_{+})\big) \cap H^{j}\big( (0,T);H^{3+4(l-j)}(\mathbf{R}_{+})\big) $ of 
\begin{eqnarray*}
\left\{
\begin{array}{ll}
\mbox{\mathversion{bold}$u$}_{t}= \alpha \mbox{\mathversion{bold}$u$}_{xxx}-\varepsilon \mbox{\mathversion{bold}$u$}_{xxxx}
+ {\rm A}(\mbox{\mathversion{bold}$w$}, \partial_{x})\mbox{\mathversion{bold}$u$} + \mbox{\mathversion{bold}$f$}, & x>0, t>0, \\
\mbox{\mathversion{bold}$u$}(x,0)=\mbox{\mathversion{bold}$u$}_{0}(x), & x>0,\\
\mbox{\mathversion{bold}$u$}(0,t)=\mbox{\mathversion{bold}$e$}, & t>0,\\
\mbox{\mathversion{bold}$u$}_{x}(0,t)=\mbox{\mathversion{bold}$0$}, & t>0,
\end{array}\right.
\end{eqnarray*}
through iteration. Now we need an estimate uniform in \( \varepsilon \). Via a standard energy method, we obtain
\begin{eqnarray*}
\begin{aligned}
\frac{1}{2}\frac{{\rm d}}{{\rm d}t}\| \mbox{\mathversion{bold}$u$}\| ^{2} &\leq -\alpha \mbox{\mathversion{bold}$u$}(0)\cdot \mbox{\mathversion{bold}$u$}_{xx}(0)
+ C\| \mbox{\mathversion{bold}$u$}\| ^{2}_{2} +  \varepsilon \mbox{\mathversion{bold}$u$}(0)\cdot \mbox{\mathversion{bold}$u$}_{xxx}(0)
-\delta \| \mbox{\mathversion{bold}$u$}_{xx}\|^{2}+ \| \mbox{\mathversion{bold}$f$}\| ^{2},\\
\frac{1}{2}\frac{{\rm d}}{{\rm d}t}\| \mbox{\mathversion{bold}$u$}_{x}\|^{2} &\leq \frac{\alpha }{2}|\mbox{\mathversion{bold}$u$}_{xx}(0)|^{2}
-\varepsilon \|\mbox{\mathversion{bold}$u$}_{xxx}\|^{2}
-\delta \| \mbox{\mathversion{bold}$u$}_{xx}\|^{2}
+ \varepsilon \mbox{\mathversion{bold}$u$}_{xx}(0)\cdot \mbox{\mathversion{bold}$u$}_{xxx}(0) + C\|\mbox{\mathversion{bold}$u$}_{xx}\|^{2} + \|\mbox{\mathversion{bold}$f$}\|^{2},\\
\frac{1}{2}\frac{{\rm d}}{{\rm d}t}\|\mbox{\mathversion{bold}$u$}_{xx}\|^{2}&\leq \frac{\alpha }{2}|\mbox{\mathversion{bold}$u$}_{xxx}(0)|^{2}
-\varepsilon \|\partial_{x}^{4}\mbox{\mathversion{bold}$u$}\|^{2}
-\delta \| \mbox{\mathversion{bold}$u$}_{xxx}\|^{2}
-\varepsilon \mbox{\mathversion{bold}$u$}_{xxx}(0)\cdot \partial ^{4}_{x} \mbox{\mathversion{bold}$u$}(0)
-\big( \mbox{\mathversion{bold}$u$}_{xxx}, {\rm A}_{0}(\mbox{\mathversion{bold}$w$})\mbox{\mathversion{bold}$u$}_{xx}\big) \\
 &\hspace*{2cm}- \big( \mbox{\mathversion{bold}$u$}_{xxx}, \mbox{\mathversion{bold}$f$}_{x}\big) + C\|\mbox{\mathversion{bold}$u$}\|^{2}_{2}.
\end{aligned}
\end{eqnarray*}
Using the equation, we can also get
\begin{eqnarray*}
-\varepsilon \mbox{\mathversion{bold}$u$}_{xxx}(0)\cdot \partial ^{4}_{x}\mbox{\mathversion{bold}$u$}(0)
= -\alpha |\mbox{\mathversion{bold}$u$}_{xxx}(0)|^{2}
-\mbox{\mathversion{bold}$u$}_{xxx}(0)\cdot \big( {\rm A}(\mbox{\mathversion{bold}$w$},\partial _{x})\mbox{\mathversion{bold}$u$}\big)(0)
-\mbox{\mathversion{bold}$u$}_{xxx}(0)\cdot \mbox{\mathversion{bold}$f$}(0).
\end{eqnarray*}
From the above estimate, we obtain 
\begin{eqnarray*}
\frac{1}{2}\frac{{\rm d}}{{\rm d}t}\|\mbox{\mathversion{bold}$u$}_{xx}\|^{2}&\leq -\frac{\alpha }{4}| \mbox{\mathversion{bold}$u$}_{xxx}(0)|^{2}
-\varepsilon \| \partial ^{4}_{x}\mbox{\mathversion{bold}$u$}\|^{2}
-\frac{\delta}{4}\| \mbox{\mathversion{bold}$u$}_{xxx}\|^{2} + C\big( \|\mbox{\mathversion{bold}$u$}\|^{2}_{2} + \|\mbox{\mathversion{bold}$f$}\|^{2}_{1}\big),
\end{eqnarray*}
which combined with the other two estimates yields
\begin{eqnarray*}
\sup _{0\leq t\leq T}\| \mbox{\mathversion{bold}$u$}(t)\|^{2}_{2} + \int ^{T}_{0} \big( \varepsilon \|\mbox{\mathversion{bold}$u$}_{xx}(t)\|^{2}_{2} 
+ \delta \|\mbox{\mathversion{bold}$u$}_{x}(t)\|^{2}_{2}\big) {\rm d}t
\leq C\bigg\{ \|\mbox{\mathversion{bold}$u$}_{0}\|^{2}_{2} + \int ^{T}_{0} \| \mbox{\mathversion{bold}$f$}(t)\|^{2}_{1}{\rm d}t \bigg\},
\end{eqnarray*}
where \( C\) is independent of \( \varepsilon \). Taking the \( t\) derivatives of the equation and estimating in the same way, we have 
for \( 0\leq m\leq l\), 
\begin{eqnarray*}
\begin{aligned}
&\sup _{0\leq t\leq T}\| \partial ^{m}_{t}\mbox{\mathversion{bold}$u$}(t)\|^{2}_{2+3(l-m)} 
+ \int ^{T}_{0} \| \partial ^{m}_{t}\mbox{\mathversion{bold}$u$}_{x}(t)\|^{2}_{2+3(l-m)}{\rm d}t \\
&\hspace*{1cm}\leq C\left\{ \|\mbox{\mathversion{bold}$u$}_{0}\|^{2}_{2+4l}+ \sum ^{l-1}_{j=0}\| \partial ^{j}_{t}\mbox{\mathversion{bold}$f$}(\cdot ,0)\|^{2}_{2+4(l-1-j)}
+\sum ^{l}_{j=0}\int ^{T}_{0}\| \partial ^{j}_{t}\mbox{\mathversion{bold}$f$}(t)\|^{2}_{1}{\rm d}t \right\}.
\end{aligned}
\end{eqnarray*}
After passing to the limit \( \varepsilon \rightarrow 0\), we obtain the solution of the limit system. Similarly as above, we see that the solution satisfies for
\( 0\leq m\leq l\),
\begin{eqnarray*}
\begin{aligned}
&\sup _{0\leq t\leq T}\| \partial_{t} ^{m}\mbox{\mathversion{bold}$u$}(t)\|^{2}_{2+3(l-m)} + \int ^{T}_{0} 
\| \partial_{t} ^{m}\mbox{\mathversion{bold}$u$}_{x}(t)\|^{2}_{2+3(l-m)}{\rm d}t \\
&\hspace*{1cm}\leq C\left\{ \|\mbox{\mathversion{bold}$u$}_{0}\|^{2}_{2+3l}+ \sum ^{l-1}_{j=0}\| \partial ^{j}_{t}\mbox{\mathversion{bold}$f$}(\cdot ,0)\|^{2}_{2+3(l-1-j)}
+\sum ^{l}_{j=0}\int ^{T}_{0}\| \partial ^{j}_{t}\mbox{\mathversion{bold}$f$}(t)\|^{2}_{1}{\rm d}t \right\}.
\end{aligned}
\end{eqnarray*}
Thus, we have obtained our second main theorem.
\begin{Th}
For any \( T>0\) and an arbitrary non-negative integer \( l\), if \( {\mbox{\mathversion{bold}$u$}}_{0}\in H^{2+3l}(\mathbf{R}_{+})\), 
\( {\mbox{\mathversion{bold}$f$}}\in Y^{l}\), and \( {\mbox{\mathversion{bold}$w$}}\in Z^{l}\) satisfy the compatibility conditions up to order \( l\), 
a unique solution \( {\mbox{\mathversion{bold}$u$}}\) of {\rm (\ref{lpd2})} exists such that \( {\mbox{\mathversion{bold}$u$}}\in X^{l}\).
\end{Th}

\end{document}